\pdfoutput=1
\RequirePackage{ifpdf}
\ifpdf
\documentclass[pdftex]{sigma}
\else
\documentclass{sigma}
\fi

\usepackage{bm,amsxtra}
\numberwithin{equation}{section}
\numberwithin{theorem}{section}
\numberwithin{proposition}{section}
\numberwithin{lemma}{section}
\numberwithin{corollary}{section}
\numberwithin{definition}{section}
\numberwithin{example}{section}
\numberwithin{remark}{section}
\numberwithin{note}{section}
\newtheorem{theoremI}{Theorem}[]

\newcommand{\Hom}{\operatorname{Hom}}
\newcommand{\Ker}{\operatorname{Ker}}
\newcommand{\ad}{\operatorname{ad}}
\newcommand{\gl}{\operatorname{\mathfrak{gl}}}

\newcommand{\csp}{\operatorname{\mathfrak{csp}}}

\newcommand{\Der}{\operatorname{Der}}
\newcommand{\Aut}{\operatorname{Aut}}
\newcommand{\gla}[1]{\mathfrak#1=\bigoplus\limits_{p\in\mathbb Z}\mathfrak#1_p} \renewcommand{\labelenumi}{{\rm (\arabic{enumi})}}

\begin{document}

\allowdisplaybreaks

\renewcommand{\PaperNumber}{038}

\FirstPageHeading

\ShortArticleName{On Free Pseudo-Product Fundamental Graded Lie Algebras}

\ArticleName{On Free Pseudo-Product\\ Fundamental Graded Lie Algebras}

\Author{Tomoaki YATSUI}

\AuthorNameForHeading{T.~Yatsui}

\Address{Department of Mathematics, Asahikawa Medical University, Asahikawa 078-8510, Japan}
\Email{\href{mailto:yatsui@asahikawa-med.ac.jp}{yatsui@asahikawa-med.ac.jp}}

\ArticleDates{Received December 16, 2011, in f\/inal form June 14, 2012; Published online June 27, 2012}

\Abstract{In this paper we f\/irst state the classif\/ication of the prolongations of complex free fundamental graded Lie algebras.
Next we introduce the notion of free pseudo-product fundamental graded Lie
algebras and study the prolongations of complex free pseudo-product fundamental graded Lie
algebras.
Furthermore we investigate the automorphism group of the prolongation of
complex free pseudo-product fundamental graded Lie algebras.}

\Keywords{fundamental graded Lie algebra; prolongation; pseudo-product graded Lie algebra}

\Classification{17B70}

\section{Introduction}

Let $\mathfrak m=\bigoplus\limits_{p<0}\mathfrak g_p$ be a graded Lie
algebra over the f\/ield $\mathbb R$ of real numbers or
the f\/ield $\mathbb C$ of complex numbers, and let $\mu$ be a positive integer.
The graded Lie algebra $\mathfrak m=\bigoplus\limits_{p<0}\mathfrak g_p $ is
called a~fundamental graded Lie algebra if the following conditions hold:
$(i)$~$\mathfrak m$ is f\/inite-dimensional; $(ii)$~$\mathfrak g_{-1}\ne\{0\}$, and
$\mathfrak m$ is generated by $\mathfrak g_{-1}$.
Moreover a fundamental graded Lie algebra
$\mathfrak m=\bigoplus\limits_{p<0}\mathfrak g_p$ is said to be
of the $\mu$-th kind if $\mathfrak g_{-\mu}\ne\{0\}$, and $\mathfrak g_p=\{0\}$ for all $p<-\mu$.
It is shown that every fundamental graded algebra $\mathfrak m=\bigoplus\limits_{p<0}\mathfrak g_p$
is prolonged to a graded Lie algebra $\gla{g (\mathfrak m)}$ satisfying the following conditions:
$(i)$~$\mathfrak g(\mathfrak m)_p=\mathfrak g_p$ for all $p<0$;
$(ii)$~for $X\in\mathfrak g(\mathfrak m)_p$ $(p\geqq0)$,
$[X,\mathfrak m]=\{0\}$ implies $X=0$;
$(iii)$~$\mathfrak g (\mathfrak m)$ is maximum among graded Lie algebras
satisfying conditions $(i)$ and $(ii)$ above.
The graded Lie algebra $\mathfrak g (\mathfrak m)$ is called the prolongation of
$\mathfrak m$.
Note that  $\mathfrak g(\mathfrak m)_0$  is the Lie algebra of all the
derivations of $\mathfrak m$ as a graded Lie algebra.

Let $\mathfrak m=\bigoplus\limits_{p<0}\mathfrak g_p$ be a fundamental graded Lie algebra of the $\mu$-th kind, where $\mu\geqq2$.
The fundamental graded Lie algebra
$\mathfrak m$ is called a free fundamental graded Lie algebra of type
$(n,\mu)$ if the following universal properties hold:
\begin{enumerate}\itemsep=0pt
\renewcommand{\labelenumi}{$(\roman{enumi})$}
\item $\dim \mathfrak g_{-1}=n$;
\item Let $\mathfrak m'=\bigoplus\limits_{p<0}\mathfrak g'_p$ be a fundamental
graded Lie algebra of the $\mu$-th kind and let $\varphi$ be a~surjective
linear mapping of $\mathfrak g_{-1}$ onto $\mathfrak g'_{-1}$.
Then $\varphi$ can be extended uniquely to a graded Lie algebra
epimorphism of $\mathfrak m$ onto $\mathfrak m'$.
\end{enumerate}
In Section~\ref{section3} we see that
a universal fundamental graded Lie algebra $b(V,\mu)$ of the $\mu$-th kind introduced
by N.~Tanaka~\cite{tan70:1} becomes
a free fundamental graded Lie algebra of type
$(n,\mu)$, where $\mu\geqq2$, and $V$ is a vector space such that
$\dim V=n\geqq2$.

In \cite{war07:1}, B.~Warhurst gave the complete list of the prolongations of
real free fundamental graded Lie algebras by using a Hall basis of a free Lie algebra.
The complex version of his theorem has the completely same form
except for the ground number f\/ield as follows:
\begin{theoremI}\label{Yatsui-theoremI}
Let $\mathfrak m=\bigoplus\limits_{p<0}\mathfrak g_p$ be a free
fundamental graded Lie algebra of type $(n,\mu)$ over $\mathbb C$.
Then
the prolongation $\gla{g (\mathfrak m)}$ of $\mathfrak m$ is
one of the following types:
\begin{enumerate}\itemsep=0pt
\renewcommand{\labelenumi}{$(\alph{enumi})$}
\item $(n,\mu)\ne(n,2)$ $(n\geqq2)$, $(2,3)$. In this case,
$\mathfrak g(\mathfrak m)_1=\{0\}$.
\item
$(n,\mu)=
(n,2)$ $(n\geqq3)$, $(2,3)$.
In this case, $\dim\mathfrak g(\mathfrak m)<\infty$ and
$\mathfrak g(\mathfrak m)_1\ne\{0\}$. Furthermore
$\mathfrak g(\mathfrak m)$ is isomorphic to a finite-dimensional simple graded
Lie algebra of type $(B_n,\{\alpha_n\})$ $(n\geqq3)$ or $(G_2,\{\alpha_1
\})$ $(n=2)$ $($see {\rm \cite{yam93:1}} or Section~{\rm \ref{section5}} for the gradations of
finite-dimensional simple graded Lie algebras over~$\mathbb C)$.
\item $(n,\mu)=(2,2)$.
In this case,
$\dim\mathfrak g(\mathfrak m)=\infty$. Furthermore,
$\mathfrak g(\mathfrak m)$ is isomorphic to the contact algebra
$K(1)$ as a graded Lie algebra.
\end{enumerate}
\end{theoremI}
The f\/irst purpose of this paper is to give a proof of Theorem~\ref{Yatsui-theoremI} by using the classif\/ication
of complex irreducible transitive graded Lie algebras of f\/inite depth
(cf.~\cite{mt70:1}).
Note that Warhurst's methods in \cite{war07:1} are available to the proof of Theorem~\ref{Yatsui-theoremI}.

Next we introduce the notion of free pseudo-product fundamental graded Lie
algebras.
Let $\mathfrak m=\bigoplus\limits_{p<0}\mathfrak g_p$ be a fundamental graded Lie algebra, and let
$\mathfrak e$ and $\mathfrak f$ be nonzero subspaces of~$\mathfrak g_{-1}$.
Then $\mathfrak m$ is called a pseudo-product fundamental graded Lie algebra
with pseudo-product structure $(\mathfrak e,\mathfrak f)$
if the following conditions hold:
$(i)$ $\mathfrak g_{-1}=\mathfrak e\oplus\mathfrak f$;
$(ii)$ $[\mathfrak e,\mathfrak e]=[\mathfrak f,\mathfrak f]=\{0\}$ (cf.~\cite{tan85:01}).

Let $\mathfrak m=\bigoplus\limits_{p<0}\mathfrak g_p$ be a pseudo-product
fundamental graded Lie algebra with a pseudo-product structure
$(\mathfrak e,\mathfrak f)$, and let
$\mathfrak g(\mathfrak m)=\bigoplus\limits_{p\in\mathbb Z}\mathfrak g(\mathfrak m)_p$ be the prolongation of $\mathfrak m$. Moreover let $\mathfrak g_0$  be the Lie algebra of all the derivations of
$\mathfrak m$ as a graded Lie algebra
preserving $\mathfrak e$ and $\mathfrak f$.
Also for $p\geqq1$ we set $\mathfrak g_p=\{X\in\mathfrak g(\mathfrak m)_p:
[X,\mathfrak g_k]\subset \mathfrak g_{p+k}\ \text{for all}\ k<0\}$ inductively.
Then the direct sum $\gla g$ becomes a graded subalgebra of $\mathfrak g(\mathfrak m)$, which is called the prolongation of
$(\mathfrak m;\mathfrak e,\mathfrak f)$.

Let $\mathfrak m=\bigoplus\limits_{p<0}\mathfrak g_p$ be a pseudo-product
fundamental graded Lie algebra of the $\mu$-th kind with pseudo-product
structure $(\mathfrak e,\mathfrak f)$, where $\mu\geqq2$.
The pseudo-product fundamental graded Lie algebra
$\mathfrak m=\bigoplus\limits_{p<0}\mathfrak g_p$ is called a free pseudo-product fundamental graded Lie algebra of type $(m,n,\mu)$
if the following conditions hold:
\begin{enumerate}\itemsep=0pt
\renewcommand{\labelenumi}{$(\roman{enumi})$}
\item $\dim\mathfrak e=m$ and $\dim\mathfrak f=n$;
\item
Let $\mathfrak m'=\bigoplus\limits_{p<0}\mathfrak g'_p$ be a pseudo-product
fundamental graded Lie algebra of the $\mu$-th kind with pseudo-product
structure $(\mathfrak e',\mathfrak f')$ and let $\varphi$ be a surjective
linear mapping of~$\mathfrak g_{-1}$ onto~$\mathfrak g'_{-1}$ such that
$\varphi(\mathfrak e)\subset\mathfrak e'$ and
$\varphi(\mathfrak f)\subset\mathfrak f'$.
Then $\varphi$ can be extended uniquely to a graded Lie algebra
epimorphism of $\mathfrak m$ onto $\mathfrak m'$.
\end{enumerate}

The main purpose of this paper is to prove the following theorem.

\begin{theoremI}\label{Yatsui-theoremII}
Let $\mathfrak m=\bigoplus\limits_{p<0}\mathfrak g_p$ be a free
pseudo-product fundamental graded Lie algebra of type $(m,n,\mu)$ with
pseudo-product structure $(\mathfrak e,\mathfrak f)$ over $\mathbb C$, and let
$\gla g$
be the prolongation of $(\mathfrak m;\mathfrak e,\mathfrak f)$.
If $\mathfrak g_{1}\ne\{0\}$, then
$\gla g$
is a finite-dimensional simple graded Lie algebra of type
$(A_{m+n},\{\alpha_m,\alpha_{m+1}\})$.
\end{theoremI}

Let $\gla g$ be the prolongation of a free pseudo-product fundamental
graded Lie algebra $\mathfrak m=\bigoplus\limits_{p<0}\mathfrak g_p$ with pseudo-product
structure $(\mathfrak e,\mathfrak f)$ over $\mathbb C$.
We denote by $\Aut(\mathfrak g;\mathfrak e,\mathfrak f)_0$ the group of
all the automorphisms as a graded Lie algebra preserving $\mathfrak e$ and
$\mathfrak f$, which is called the automorphism group of
the pseudo-product graded Lie algebra $\gla g$.
In Section~\ref{section9}, we show that $\Aut(\mathfrak g;\mathfrak e,\mathfrak f)_0$
is isomorphic to
$GL(\mathfrak e)\times GL(\mathfrak f)$.

\subsection*{Notation and conventions}
\begin{enumerate}\itemsep=0pt
\item From Section~\ref{section2} to the last section, all vector spaces are considered over
the f\/ield $\mathbb C$ of complex numbers.
\item
Let $V$ be a vector space and let $W_1$ and $W_2$ be subspaces of $V$.
We denote by $W_1\wedge W_2$ the subspace of $\Lambda^2V$
spanned by all the elements of the form $w_1\wedge w_2$
$(w_1\in W_1,w_2\in W_2)$.
\item Graded vector spaces are always $\mathbb Z$-graded. If we write
$V=\bigoplus\limits_{p<0} V_p$, then
it is understood that $V_p=\{0\}$ for all $p\geqq0$.
Let $V=\bigoplus\limits_{p\in\mathbb Z}V_p$ be a graded vector space.
We denote by $V_-$ the subspace $V=\bigoplus\limits_{p<0} V_p$. Also
for $k\in \mathbb Z$ we denote by
$V_{\leqq k}$ the subspace $\bigoplus\limits_{p\leqq k} V_p$.

Let $V=\bigoplus\limits_{p\in\mathbb Z}V_p$ and
$W=\bigoplus\limits_{p\in\mathbb Z}W_p$ be graded vector spaces.
For $r\in\mathbb Z$, we set
\[
\Hom(V,W)_r=\{\varphi\in\Hom(V,W):\varphi(V_p)\subset W_{p+r}\
\text{for all} \ p\in\mathbb Z\}.
\]
\end{enumerate}

\section{Free fundamental graded Lie algebras}\label{section2}

First of all we give several def\/initions about graded Lie algebras.
Let $\mathfrak g$ be a Lie algebra. Assume that there is given a family of
subspaces $(\mathfrak g_p)_{p\in\mathbb Z}$ of $\mathfrak g$ satisfying
the following conditions:
\begin{enumerate}\itemsep=0pt
\renewcommand{\labelenumi}{$(\roman{enumi})$}
\item $\gla g$;
\item $\dim \mathfrak g_p<\infty$ for all $p\in\mathbb Z$;
\item $[\mathfrak g_p,\mathfrak g_q]\subset \mathfrak g_{p+q}$
for all $p,q\in\mathbb Z$.
\end{enumerate}
Under these conditions, we say that $\gla g$ is a graded Lie algebra (GLA).
Moreover we def\/ine the notion of homomorphism, isomorphism, monomorphism,
epimorphism, subalgebra and
ideal for GLAs in an obvious manner.

A GLA $\gla g$ is called transitive if
for $X\in\mathfrak g_p$ $(p\geqq0)$, $[X,\mathfrak g_-]=\{0\}$ implies $X=0$,
where $\mathfrak g_-$ is the negative part
$\bigoplus\limits_{p<0}\mathfrak g_p$ of $\mathfrak g$.
Furthermore a GLA $\gla g$ is called irreducible if
the $\mathfrak g_0$-module $\mathfrak g_{-1}$ is irreducible.

Let $\mu$ be a positive integer. A GLA
$\gla g$ is said to be of depth $\mu$ if $\mathfrak g_{-\mu}\ne\{0\}$ and
$\mathfrak g_p=\{0\}$ for all $p<-\mu$.

Next we def\/ine fundamental GLAs.
A GLA $\mathfrak m=\bigoplus\limits_{p<0}\mathfrak g_p$ is
called a \textit{fundamental graded Lie algebra} (FGLA)
if the following conditions hold:
\begin{enumerate}\itemsep=0pt
\renewcommand{\labelenumi}{$(\roman{enumi})$}
\item $\dim \mathfrak m<\infty$;
\item $\mathfrak g_{-1}\ne\{0\}$, and $\mathfrak m$ is generated by $\mathfrak g_{-1}$,
or more precisely $\mathfrak g_{p-1}=[\mathfrak g_p,\mathfrak g_{-1}]$
for all $p<0$.
\end{enumerate}
If an FGLA $\mathfrak m=\bigoplus\limits_{p<0}\mathfrak g_p$
is of depth $\mu$, then $\mathfrak m$
is also said to be of the $\mu$-th kind.
Moreover an FGLA $\mathfrak m=\bigoplus\limits_{p<0}\mathfrak g_p$ is called non-degenerate if
for $X\in\mathfrak g_{-1}$, $[X,\mathfrak g_{-1}]=\{0\}$ implies $X=0$.
\par
Let $\mathfrak m=\bigoplus\limits_{p<0}\mathfrak g_p$ be an FGLA of
the $\mu$-th kind, where $\mu\geqq2$. $\mathfrak m$ is called a free
fundamental graded Lie algebra of type $(n,\mu)$
if the following conditions hold:
\begin{enumerate}\itemsep=0pt
\renewcommand{\labelenumi}{$(\roman{enumi})$}
\item $\dim \mathfrak g_{-1}=n$;
\item Let $\mathfrak m'=\bigoplus\limits_{p<0}\mathfrak g'_p$ be an FGLA
of the $\mu$-th kind and let $\varphi$ be a surjective linear
mapping of $\mathfrak g_{-1}$ onto $\mathfrak g'_{-1}$.
Then $\varphi$ can be extended uniquely to a GLA
epimorphism of $\mathfrak m$ onto $\mathfrak m'$.
\end{enumerate}

\begin{proposition}\label{prop2.1} Let $n$ and $\mu$ be positive integers such that $n,\mu\geqq2$.
\begin{enumerate}\itemsep=0pt
\item There exists a unique free FGLA of type
$(n,\mu)$ up to isomorphism.
\item Let $\mathfrak m=\bigoplus\limits_{p<0}\mathfrak g_p$ be a free
FGLA of type $(n,\mu)$.
We denote by $\Der(\mathfrak m)_0$ the Lie algebra of all the derivations
of $\mathfrak m$ preserving the gradation of $\mathfrak m$.
Then the mapping $\Phi:\Der(\mathfrak m)_0\ni D\mapsto D|\mathfrak g_{-1}\in\gl(\mathfrak g_{-1})$ is a Lie algebra isomorphism.
\end{enumerate}
\end{proposition}
\begin{proof} (1)
The uniqueness of a free FGLA of type $(n,\mu)$ follows from the def\/inition.
We set $X=\{1,\dots,n\}$. Let $L(X)$ be the free Lie algebra on $X$
(see \cite[Chapter~II, \S~2]{bou72:1})
and let $i:X\to L(X)$ be the canonical injection.
We def\/ine a mapping $\phi$ of $X$ into $\mathbb Z$ by $\phi(k)=-1$ $(k\in X)$.
The mapping $\phi$ def\/ines the natural gradation $(L(X)_p)_{p<0}$
on $L(X)$ such that: $(i)$ $L(X)$ is generated by $L(X)_{-1}$;
$(ii)$ $\{i(1),\dots,i(n)\}$ is a basis of $L(X)_{-1}$
(see \cite[Chapter~II, \S~2, no.~6]{bou72:1}). Note that if $n>1$, then
$L(X)_p\ne0$ for all $p<0$.
We set $\mathfrak a=\bigoplus\limits_{p<-\mu}L(X)_p$;
then $\mathfrak a$ is a graded ideal of $L(X)$ and the factor GLA $\mathfrak m=L(X)/\mathfrak a$ becomes an FGLA of the $\mu$-th kind.
We put $\mathfrak a_p=\mathfrak a\cap L(X)_p$ and
$\mathfrak g_p=L(X)_p/\mathfrak a_p$.

Now we prove that $\mathfrak m=\bigoplus\limits_{p<0}\mathfrak g_p$ is a free
FGLA of type $(n,\mu)$.
Let $\mathfrak m'=\bigoplus\limits_{p<0}\mathfrak g'_p$ be an FGLA of the
$\mu$-th kind and let
$\varphi$ be a surjective linear mapping of $\mathfrak g_{-1}$ onto
$\mathfrak g'_{-1}$.
Let $h$ be a mapping of~$X$ into $\mathfrak m'$ def\/ined by $h(k)=\varphi(i(k))$ $(k\in X)$.
Then there exists a Lie algebra homomorphism~$\tilde{h}$ of~$L(X)$
into $\mathfrak m'$ such that $\tilde{h}\circ i=h$.
Since $L(X)$ (resp.~$\mathfrak m'$) is generated by
$L(X)_{-1}$ (resp.~$\mathfrak g'_{-1}$), $\tilde{h}$ is surjective.
Since $\mathfrak m'=\bigoplus\limits_{p<0}\mathfrak g'_p$ is of the
$\mu$-th kind, $\tilde{h}(\mathfrak a)=0$, so
$\tilde{h}$ induces a GLA epimorphism~$L(\varphi)$ of
$\mathfrak m$ onto~$\mathfrak m'$ such that
$L(\varphi)|\mathfrak g_{-1}=\varphi$.
The homomorphism~$L(\varphi)$ is unique, because
$\mathfrak m=\bigoplus\limits_{p<0}\mathfrak g_p$ is generated by
$\mathfrak g_{-1}$.
Thus $\mathfrak m$ is a free FGLA of type~$(n,\mu)$.

(2) Assume that $\mathfrak m$ is a free FGLA
constructed in~(1).
Let $\phi$ be an endomorphism of $\mathfrak g_{-1}$.
By Corollary to Proposition~8 of \cite[Chapter~II, \S~2, no.~8]{bou72:1},
$\phi$ can be extended uniquely to a~unique derivation $D$ of $L(X)$.
Since $D(L(X)_{-1})=\phi(L(X)_{-1})=\phi(\mathfrak g_{-1})\subset L(X)_{-1}$,
and since~$L(X)$ is generated by $L(X)_{-1}$, we see that
$D(L(X)_{p})\subset L(X)_{p}$ and $D(\mathfrak a)\subset \mathfrak a$.
Thus there is a~derivation of $D_\phi$ of $\mathfrak m$ such that
$\pi\circ D=D_\phi\circ \pi$, where $\pi$ is the natural projection of $L(X)$ onto $\mathfrak m$.
The correspondence
$\gl(\mathfrak g_{-1})\ni\phi\mapsto D_\phi\in\Der(\mathfrak m)_0$
is an injective linear mapping.
Hence $\dim \gl(\mathfrak g_{-1})\leqq\dim \Der(\mathfrak m)_0$.
On the other hand, since $\mathfrak m$ is generated by $\mathfrak g_{-1}$,
the mapping $\Phi$ is a~Lie algebra monomorphism.
Therefore $\Phi$ is a Lie algebra isomorphism.
\end{proof}

\begin{remark}
Let $n$ and $\mu$ be positive integers with $n,\mu\geqq2$, and
let $\mathfrak m=\bigoplus\limits_{p<0}\mathfrak g_p$ be a free FGLA of
type $(n,\mu)$.  Furthermore
let $\mathfrak m'=\bigoplus\limits_{p<0}\mathfrak g'_p$
be an FGLA of the $\mu$-th kind,
and let $\varphi$ be a linear mapping of $\mathfrak g_{-1}$ into
$\mathfrak g'_{-1}$.
\begin{enumerate}\itemsep=0pt
\item
From the proof of Proposition \ref{prop2.1}, there exists a unique GLA
homomorphism $L(\varphi)$ of $\mathfrak m$ into $\mathfrak m'$
such that $L(\varphi)|\mathfrak g_{-1}=\varphi$.
\item
Let $\mathfrak m''=\bigoplus\limits_{p<0}\mathfrak g''_p$
be an FGLA of the $\mu$-th kind,
and let $\varphi'$ be a linear mapping of $\mathfrak g'_{-1}$ into~$\mathfrak g''_{-1}$. Assume that $\mathfrak m'=\bigoplus\limits_{p<0}\mathfrak g'_p$ is a free FGLA. By the uniqueness of $L(\varphi'\circ\varphi)$,
we see that
$L(\varphi'\circ\varphi)=L(\varphi')\circ L(\varphi)$.
\item Assume that $\mathfrak m'=\bigoplus\limits_{p<0}\mathfrak g'_p$ is a
free FGLA and $\varphi$ is injective. By the result of (2),
$L(\varphi)$ is a~monomorphism.
\item Let $W$ be an $m$-dimensional subspace of $\mathfrak g_{-1}$ with $m\geqq2$.
By the result of (3), the subalgebra of $\mathfrak m$ generated by $W$ is a
free FGLA of type $(m,\mu)$.
\end{enumerate}
\end{remark}
By Remark 2.1 (4) and \cite[Chapter~II, \S~2, Theorem~1]{bou72:1}, we get the following lemma.
\begin{lemma}\label{lem2.1} Let $\mathfrak m=\bigoplus\limits_{p<0}\mathfrak g_p$ be a
free FGLA of type $(n,\mu)$ with $\mu\geqq3$.
If $X$, $Y$ are linearly independent elements of $\mathfrak g_{-1}$, then
\begin{gather*}
\ad(X)^\mu(Y)=0, \qquad \ad(X)^{\mu-1}(Y)\ne 0, \\
\ad(Y)\ad(X)^{\mu-1}(Y)=0, \qquad \ad(Y)\ad(X)^{\mu-2}(Y)\ne 0.
\end{gather*}
\end{lemma}

\section{Universal fundamental graded Lie algebras}\label{section3}

Following N.~Tanaka \cite{tan70:1}, we introduce universal FGLAs of the
$\mu$-th kind.

Let $V$ be an $n$-dimensional vector space.
We def\/ine vector spaces $b(V)_p$ $(p<0)$ and li\-near mappings $B_p$ of
$\sum\limits_{r+s=p}b(V)_r\wedge b(V)_s$ into $b(V)_p$ $(p\leqq -2)$
as follows:
First of all, we put \mbox{$b(V)_{-1}=V$} and $b(V)_{-2}=\Lambda^2V$. Further
we def\/ine a mapping $B_{-2}:b(V)_{-1}\wedge b(V)_{-1}\to b(V)_{-2}$
to be the identity mapping. For $k\leqq-3$, we def\/ine $b(V)_k$ and
$B_k$ inductively as follows:
We set $b(V)^{(k+1)}=\bigoplus\limits_{p=-1}^{k+1}b(V)_p$ and
we def\/ine a subspace $c(V)_k$ of $\Lambda^2(b(V)^{(k+1)})$ to be
$\sum\limits_{r+s=k}b(V)_r\wedge b(V)_s$.
We denote by $A(V)_k$ the subspace of $c(V)_k$ spanned by the elements
\[
\underset{(X,Y,Z)}{\mathfrak S}\sum_{r+s=k}\sum_{u+v=r}B_r(X_u\wedge Y_v)\wedge Z_s,
\qquad  X,Y,Z\in b(V)^{(k+1)} ,
\]
where $\underset{(X,Y,Z)}{\mathfrak S}$ stands for the cyclic sum with respect to $X$, $Y$, $Z$, and $X_u$ denotes the $b(V)_u$-component in the decomposition
$b(V)^{(k+1)}=\bigoplus\limits_{p=-1}^{k+1}b(V)_p$.
Now we def\/ine $b(V)_k$ to be the factor space $c(V)_k/A(V)_k$, and
$B_k$ to be the projection of $c(V)_k$ onto $b(V)_k$.
We put $b(V)=\bigoplus\limits_{p<0}b(V)_p$ and def\/ine a bracket operation
$[\ , \ ]$ on $b(V)$ by
\[
[X,Y]=\sum_{p\leqq-2}\sum_{r+s=p}B_p(X_r\wedge Y_s)
\]
for all $X,Y\in b(V)$. Then $b(V)=\bigoplus\limits_{p<0}b(V)_p$ becomes a GLA generated by $b(V)_{-1}$, and $b(V)_p\ne0$ for all $p<0$ if $\dim V>1$.

Note that $b(V)_{-3}$ is isomorphic to $\Lambda^2(V)\otimes V/\Lambda^3V$.
Let $\mu$ be a positive integer. Assume that $\mu\geqq 2$ and $\dim V=n\geqq 2$.
Since $\bigoplus\limits_{p<-\mu}b(V)_p$ is a graded ideal of $b(V)$,
we see that the factor space $b(V,\mu)=b(V)/\bigoplus\limits_{p<-\mu}b(V)_p$
becomes an FGLA of $\mu$-th kind,
which is called a universal fundamental graded Lie algebra of the
$\mu$-th kind.
By \cite[Proposition~3.2]{tan70:1},
$b(V,\mu)$ is a free FGLA of type $(n,\mu)$.

\section{The prolongations of fundamental graded Lie algebras}\label{section4}

Following N.~Tanaka \cite{tan70:1}, we introduce the prolongations of FGLAs.
Let $\mathfrak m=\bigoplus\limits_{p<0}\mathfrak g_p$ be an FGLA.
A GLA
$\mathfrak g(\mathfrak m)
=\bigoplus\limits_{p\in\mathbb Z}\mathfrak g(\mathfrak m)_p$ is called the prolongation of $\mathfrak m$ if the following conditions hold:
\begin{enumerate}\itemsep=0pt
\renewcommand{\labelenumi}{$(\roman{enumi})$}
\item $\mathfrak g(\mathfrak m)_p=\mathfrak g_p$ for all $p<0$;
\item $\mathfrak g(\mathfrak m)$ is a transitive GLA;
\item If $\mathfrak h=\bigoplus\limits_{p\in\mathbb Z}\mathfrak h_p$ is a
GLA satisfying conditions (i) and (ii) above, then
$\mathfrak h=\bigoplus\limits_{p\in\mathbb Z}\mathfrak h_p$ can be embedded in
$\mathfrak g(\mathfrak m)$ as a GLA.
\end{enumerate}
We construct the prolongation $\mathfrak g(\mathfrak m)
=\bigoplus\limits_{p\in\mathbb Z}\mathfrak g(\mathfrak m)_p$ of $\mathfrak m$.
We set $\mathfrak g(\mathfrak m)_p=\mathfrak g_p$ $(p<0)$.
We def\/ine subspaces $\mathfrak g(\mathfrak m)_k$ $(k\geqq0)$ of
$\Hom(\mathfrak m,\bigoplus\limits_{p\leqq k-1}\mathfrak g(\mathfrak m)_p)_k$
and a bracket operation on $\mathfrak g(\mathfrak m)
=\bigoplus\limits_{p\in\mathbb Z}\mathfrak g(\mathfrak m)_p$ inductively.
First $\mathfrak g(\mathfrak m)_0$ is def\/ined to be $\Der(\mathfrak m)_0$ and a bracket operation $[\ ,\ ]:\bigoplus\limits_{p\leqq 0}\mathfrak g(\mathfrak m)_p\times \bigoplus\limits_{p\leqq 0}\mathfrak g(\mathfrak m)_p\to
\bigoplus\limits_{p\leqq 0}\mathfrak g(\mathfrak m)_p$ is def\/ined by
\begin{gather*}
[X,Y]=-[Y,X]=X(Y),\qquad  X\in\mathfrak g(\mathfrak m)_0, \quad Y\in\mathfrak m,\\
[X,Y]=XY-YX, \qquad  X,Y\in\mathfrak g(\mathfrak m)_0 .
\end{gather*}
Next for $k>0$ we def\/ine $\mathfrak g(\mathfrak m)_k$ $(k\geqq1)$ inductively
as follows:
\begin{gather*}
\mathfrak g(\mathfrak m)_k
=\Big\{ X\in\Hom\Big(\mathfrak m,\!\bigoplus_{p\leqq k-1}\mathfrak g(\mathfrak m)_p\Big)_k\!:
X([u,v])=[X(u),v]+[u,X(v)]\  \text{for all}\ u,v\in\mathfrak m\Big\},
\end{gather*}
where for $X\in\mathfrak g(\mathfrak m)_r$, $u\in\mathfrak m$,
we set $[X,u]=-[u,X]=X(u)$. Further for $X\in\mathfrak g(\mathfrak m)_k$,
$Y\in\mathfrak g(\mathfrak m)_l$ $(k,l\geqq0)$, by induction on $k+l\geqq0$,
we def\/ine $[X,Y]\in\Hom(\mathfrak m,\mathfrak g(\mathfrak m))_{k+l}$ by
\[
[X,Y](u)=[X,[Y,u]]-[Y,[X,u]],\qquad  u\in \mathfrak m .
\]
It follows easily that $[X,Y]\in \mathfrak g(\mathfrak m)_{k+l}$.
With this bracket operation, $\mathfrak g(\mathfrak m)
=\bigoplus\limits_{p\in\mathbb Z}\mathfrak g(\mathfrak m)_p$ becomes a graded
Lie algebra satisfying conditions $(i)$, $(ii)$ and $(iii)$ above.

Let $\mathfrak m$ and $\mathfrak g(\mathfrak m)$ be as above.
Assume that we are given a subalgebra $\mathfrak g_0$ of
$\mathfrak g(\mathfrak m)_0$.
We def\/ine subspaces $\mathfrak g_k$ $(k\geqq 1)$ of
$\mathfrak g(\mathfrak m)_k$ inductively as follows:
\[
\mathfrak g_k=\{ X\in\mathfrak g(\mathfrak m)_k:
[X,\mathfrak g_p]\subset \mathfrak g_{p+k}\
\text{for all}\ p<0 \}.
\]
If we put $\gla g$, then
$\gla g$ becomes a transitive graded Lie subalgebra of $\mathfrak g(\mathfrak m)$,
which is called the prolongation of
$(\mathfrak m,\mathfrak g_0)$.

By Proposition \ref{prop2.1} (2) we get the following proposition.
\begin{proposition}\label{prop4.1} Let $\mathfrak m=\bigoplus\limits_{p<0}\mathfrak g_p$ be a
free FGLA and let $\mathfrak g(\mathfrak m)
=\bigoplus\limits_{p\in\mathbb Z}\mathfrak g(\mathfrak m)_p$
be the prolongation of $\mathfrak m$.
Then the mapping $\mathfrak g(\mathfrak m)_0\ni D\mapsto D|\mathfrak g_{-1}
\in\gl(\mathfrak g_{-1})$ is an isomorphism.
\end{proposition}
Conversely we obtain the following proposition.

\begin{proposition}\label{prop4.2}
Let $\mathfrak m=\bigoplus\limits_{p<0}\mathfrak g_p$ be an
FGLA of the $\mu$-th kind and let
$\mathfrak g(\mathfrak m)
=\bigoplus\limits_{p\in\mathbb Z}\mathfrak g(\mathfrak m)_p$
be the prolongation of $\mathfrak m$.
Assume that $\mathfrak g(\mathfrak m)_0$ is isomorphic to
$\gl(\mathfrak g_{-1})$.
If $\mu=2$ or $\mu=3$,
then $\mathfrak m$ is a~free FGLA.
\end{proposition}

\begin{proof}\looseness=-1
We put $n=\dim \mathfrak g_{-1}$. We consider a universal FGLA
$b(\mathfrak g_{-1},\mu)=\bigoplus\limits_{p<0}b(\mathfrak g_{-1},\mu)_p$ of
the $\mu$-th kind. Since $b(\mathfrak g_{-1},\mu)$
is a free FGLA of type $(n,\mu)$,
there exists a GLA epimorphism $\varphi$
of $b(\mathfrak g_{-1},\mu)$ onto $\mathfrak m$ such that the restriction
$\varphi|b(\mathfrak g_{-1},\mu)_{-1}$ is the identity mapping.
Let $\check{b}(\mathfrak g_{-1},\mu)
=\bigoplus\limits_{p\in\mathbb Z}\check{b}(\mathfrak g_{-1},\mu)_p$
be the prolongation of $b(\mathfrak g_{-1},\mu)$.
Since the mapping $\mathfrak g(\mathfrak m)_0\ni D\mapsto
D|\mathfrak g_{-1}\in\gl(\mathfrak g_{-1})$ is an isomorphism,
$\varphi$ can be extended to be a homomorphism $\check{\varphi}$
of $\bigoplus\limits_{p\leqq0}\check{b}(\mathfrak g_{-1},\mu)_p$ onto
$\bigoplus\limits_{p\leqq0}\mathfrak g(\mathfrak m)_p$.
Let $\mathfrak a$ be the kernel of $\check{\varphi}$;
then $\mathfrak a$ is a graded ideal of $\bigoplus\limits_{p\leqq0}\check{b}(\mathfrak g_{-1},\mu)_p$. We set $\mathfrak a_p
=\mathfrak a\cap \check{b}(\mathfrak g_{-1},\mu)_p$; then
$\mathfrak a=\bigoplus\limits_{p\leqq 0}\mathfrak a_p$.
Since the restriction of $\check{\varphi}$ to $\check{b}(\mathfrak g_{-1},\mu)_{-1}
\oplus \check{b}(\mathfrak g_{-1},\mu)_0$ is injective,
$\mathfrak a_p=\{0\}$ for $p\geqq-1$. Also
each $\mathfrak a_p$ is a $\check{b}(\mathfrak g_{-1},\mu)_0$-submodule of
$\check{b}(\mathfrak g_{-1},\mu)_p$.
From the construction of $b(\mathfrak g_{-1},\mu)$, we see that
$b(\mathfrak g_{-1},\mu)_{-2}$ (resp. $b(\mathfrak g_{-1},\mu)_{-3}$) is
isomorphic to
$\Lambda^2(\mathfrak g_{-1})$
(resp.\ $\Lambda^2(\mathfrak g_{-1})\otimes \mathfrak g_{-1}/
\Lambda^3(\mathfrak g_{-1}))$ as a $\check{b}(\mathfrak g_{-1},\mu)_0$-module.
By the table of \cite{ov90:1},
$\Lambda^2(\mathfrak g_{-1})$ and $\Lambda^2(\mathfrak g_{-1})\otimes
\mathfrak g_{-1}/\Lambda^3(\mathfrak g_{-1})$ are irreducible
$\gl(\mathfrak g_{-1})$-modules.
Thus we see that $\mathfrak a_{-2}=\mathfrak a_{-3}=\{0\}$.
From $\mu\leqq3$ it follows that $\varphi$ is an isomorphism.
\end{proof}

\section{Finite-dimensional simple graded Lie algebras}\label{section5}

Following \cite{yam93:1}, we f\/irst state the classif\/ication of
f\/inite-dimensional simple GLAs.

\looseness=-1
Let $\gla g$ be a f\/inite-dimensional simple GLA of the $\mu$-th kind over
$\mathbb C$ such that the negative part $\mathfrak g_-$ is an FGLA.
Let $\mathfrak h$ be a Cartan subalgebra of $\mathfrak g_0$;
then $\mathfrak h$ is a Cartan subalgebra of $\mathfrak g$ such that
$E\in \mathfrak h$, where
$E$ is  the element of $\mathfrak g_0$ such that
$[E,x]=px$ for all $x\in\mathfrak g_p$ and~$p$.
Let~$\Delta$ be
a root system of $(\mathfrak g,\mathfrak h)$. For~$\alpha\in\Delta$,
we denote by $\mathfrak g^{\alpha}$ the root space corresponding to~$\alpha$. We set
$\mathfrak h_\mathbb R=\{ h\in\mathfrak h:\alpha(h)\in\mathbb R\ \text{for all}\ \alpha\in \Delta \}$
and let $(h_1,\dots,h_l)$ be a basis of $\mathfrak h_\mathbb R$ such that
$h_1=E$.
We def\/ine the set of positive roots $\Delta^+$ as the set of roots
which are positive with respect to the lexicographical ordering
in $\mathfrak h_{\mathbb R}^*$ determined by the basis $(h_1,\dots,h_l)$ of
$\mathfrak h_{\mathbb R}$.
Let $\Pi\subset \Delta^+$ be the corresponding simple root system.
We denote by $\{m_1,\dots,m_l\}$ the coordinate functions
corresponding to $\Pi$, i.e., for $\alpha\in\Delta$,
we can write $\alpha=\sum\limits_{i=1}^l m_i(\alpha)\alpha_i$.

We set $\alpha_i(E)=s_i$ and $\bm s=(s_1,\dots,s_l)$;
then each $s_i$ is a non-negative integer.
For $\alpha\in\Delta$, we call the integer
$\ell_{\bm s}(\alpha)=\sum\limits_{i=1}^l m_i(\alpha)s_i$ the
$\bm s$-length of $\alpha$.
We put $\Delta_p=\{ \alpha\in\Delta: \ell_{\bm s}(\alpha)=p \}$,
$\Pi_p=\Delta_p\cap \Pi$ and $I=\{ i\in\{1,\dots,l\}: s_i=1 \}$.
Let $\theta$ be the highest root of $\mathfrak g$; then
$\ell_{\bm s}(\theta)=\mu$. Also since the $\mathfrak g_0$-module
$\mathfrak g_{-\mu}$ is irreducible,
$\dim \mathfrak g_{-\mu}=1$
if and only if $\langle \theta,\alpha\spcheck_i\rangle=0$ for all
$i\in\{1,\dots,l\}\setminus I$, where $\{\alpha\spcheck_i\}$ is the simple root system of the dual root system $\Delta\spcheck$ of $\Delta$ corresponding
to $\{\alpha_i\}$.
In our situation, since $\mathfrak g_-$ is generated by
$\mathfrak g_{-1}$, we have $s_i=0$ or 1 for all $i$.
The $l$-tuple $\bm s=(s_1,\dots,s_l)$ of non-negative integers
is determined only by the ordering of $(\alpha_1,\dots,\alpha_l)$.
In what follows, we assume that the ordering of
$(\alpha_1,\dots,\alpha_l)$ is as in the table of~\cite{bou68:1}.
If $\mathfrak g$ has the Dynkin diagram of type $X_l$ $(X=A,\dots,G)$,
then the simple GLA $\gla g$
is said to be of type $(X_l,\Pi_1)$.
Here we remark that for an automorphism $\bar{\mu}$ of the
Dynkin diagram, a~simple GLA of type $(X_l,\Pi_1)$
is isomorphic to that of type $(X_l,\bar{\mu}(\Pi_1))$.
We will identify a~simple GLA of type $(X_l,\Pi_1)$
with that of type $(X_l,\bar{\mu}(\Pi_1))$.

For $i\in I$, we put $\Delta_p^{(i)}=\{\alpha\in\Delta:
m_i(\alpha)=p\ \text{and}\ m_j(\alpha)=
0\ \text{for}\ j\in I\setminus\{i\} \}$ and
$\mathfrak g_p^{(i)}=\sum\limits_{\alpha\in\Delta_p^{(i)}}\mathfrak g^\alpha$;
then
$\mathfrak g_{-1}^{(i)}$ is an irreducible
$\mathfrak g_0$-submodule of $\mathfrak g_{-1}$ with
highest weight $-\alpha_i$.
In particular, if the $\mathfrak g_0$-module $\mathfrak g_{-1}$
is irreducible, then $\#(I)=1$.

{\samepage For $i\in I$,
we denote by $\mathfrak g^{(i)}$ the subalgebra of $\mathfrak g$ generated by
$\mathfrak g_{-1}^{(i)}\oplus \mathfrak g_{1}^{(i)}$;
then $\mathfrak g^{(i)}$ is a~simple GLA whose
Dynkin diagram is the connected component containing the vertex $i$ of
the sub\-diagram of~$X_l$
corresponding to vertices $(\{1,\dots,l\}\setminus I)\cup\{i\}$.
We denote by $\theta^{(i)}$ the highest root of $\mathfrak g^{(i)}$.
Then $[\mathfrak g_{-1}^{(i)},\mathfrak g_{-1}^{(i)}]=\{0\}$ if and only if
$m_i(\theta^{(i)})=1$.

}

From Theorem 5.2 of \cite{yam93:1}, we obtain the following theorem:
\begin{theorem}\label{thm5.1}
Let $\gla g$ be a finite-dimensional simple GLA over $\mathbb C$ such that
$\mathfrak g_-$ is an FGLA and
the $\mathfrak g_0$-module $\mathfrak g_{-1}$
is irreducible.
Then $\gla g$ is the prolongation of $\mathfrak g_-$ except for
the following cases:
\begin{enumerate}\itemsep=0pt
\renewcommand{\labelenumi}{$(\alph{enumi})$}
\item $\mathfrak g_-$ is of the first kind;
\item $\mathfrak g_-$ is of the second kind and $\dim\mathfrak g_{-2}=1$.
\end{enumerate}
\end{theorem}

Let $\gla g$ be a f\/inite-dimensional simple GLA.
Now we assume that $\mathfrak g_0$ is isomorphic to $\gl(\mathfrak g_{-1})$;
then the $\mathfrak g_0$-module $\mathfrak g_{-1}$ is irreducible.
The derived subalgebra $[\mathfrak g_0,\mathfrak g_0]$ of $\mathfrak g_0$
is a~semisimple Lie algebra whose Dynkin diagram is the subdiagram of
$X_l$ consisting
of the vertices $\{1,\dots,l\}\setminus I$.
Since $[\mathfrak g_0,\mathfrak g_0]$ is of type $A_{l-1}$ and since
the $\mathfrak g_0$-module $\mathfrak g_{-1}$ is elementary,
$(X_l,\Delta_1)$ is one of the following cases:
\[
(A_l,\{\alpha_1\}),\qquad (B_l,\{\alpha_l\}), \quad  l\geqq2,
\qquad (G_2,\{\alpha_1\}).
\]
From this result and Propositions \ref{prop4.1} and \ref{prop4.2}, we get the following theorem:
\begin{theorem}\label{thm5.2}
Let $\gla g$ be a finite-dimensional simple GLA of type
$(X_l,\Pi_1)$ over $\mathbb C$ satisfying the following conditions:
\begin{enumerate}\itemsep=0pt
\renewcommand{\labelenumi}{$(\roman{enumi})$}
\item $\mathfrak g_-$ is an FGLA of the $\mu$-th kind;
\item The $\mathfrak g_0$-module $\mathfrak g_{-1}$ is irreducible;
\item $\mathfrak g_0$ is isomorphic to $\gl(\mathfrak g_{-1})$;
\item $\mathfrak g$ is the prolongation of $\mathfrak g_{-}$.
\end{enumerate}

Then $\mathfrak g_-$ is a free FGLA of type
$(l,\mu)$, and
$\gla g$ is one of the following types:
\begin{enumerate}\itemsep=0pt
\renewcommand{\labelenumi}{$(\alph{enumi})$}
\item $l\geqq3$, $\mu=2$, $(X_l,\Pi_1)=(B_l,\{\alpha_l\})$.
\item $l=2$, $\mu=3$, $(X_l,\Pi_1)=(G_2,\{\alpha_1\})$.
\end{enumerate}
\end{theorem}

\section[Graded Lie algebras $W(n)$, $K(n)$ of Cartan type]{Graded Lie algebras $\boldsymbol{W(n)}$, $\boldsymbol{K(n)}$ of Cartan type}\label{section6}

In this section, following V.G.~Kac \cite{kac68:1},
we describe Lie algebras $W(n)$, $K(n)$ of Cartan type and their standard gradations.

Let $A(m)$ denote the monoid (under addition) of all $m$-tuples of
non-negative integers.
For an $m$-tuple
$\bm s=(s_1,\dots,s_m)$ of positive integers and
$\alpha=(\alpha_1,\dots,\alpha_m)\in A(m)$
we set $\|\alpha\|_{\bm s}=\sum\limits_{i=1}^ms_i\alpha_i$.
Also we denote the $m$-tuple $(1,\dots,1)$ by ${\bm 1}_m$
and  we denote the
$(m+1)$-tuple $(1,\dots,1,2)$ by $({\bm 1}_m,2)$.
Let $\mathfrak A(m)=\mathbb C[x_1,\dots,x_m]$.
For any $m$-tuple $\bm s$
of positive integers, we denote by $\mathfrak A(m;\bm s)_p$
the subspace of $\mathfrak A(m)$ spanned by polynomials
\[
x^{\alpha}=x_1^{\alpha_1}\cdots x_m^{\alpha_m}, \qquad
 \alpha=(\alpha_1,\dots,\alpha_m)\in A(m),\quad  \|\alpha\|_{\bm s}=p .
\]
Let $W(m)$ be the Lie algebra consisting of all the polynomial vector f\/ields
\begin{gather}\label{eq6.1}
\sum\limits_{i=1}^m P_i\frac{\partial}{\partial x_i},
\qquad  P_i\in\mathfrak A(m) .
\end{gather}

For an $m$-tuple $\bm s=(s_1,\dots,s_m)$ of positive integers,
we denote by $W(m;\bm s)_p$
the subspaces of $W(m)$ consisting of those polynomial vector f\/ields \eqref{eq6.1}
such that the polynomials $P_i$ are contained in
$\mathfrak A(m;\bm s)_{p+s_i}$;
then
$W(m;\bm s)=\bigoplus\limits_{p\in\mathbb Z}W(m;\bm s)_p$
is a transitive GLA.
In particular, $W(m;{\bm 1}_m)=\bigoplus\limits_{p\geqq-1}W(m;{\bm 1}_m)_p$ is
a transitive irreducible GLA
such that: $(i)$ $W(m;{\bm 1}_m)_0$ is isomorphic to $\gl(m,\mathbb C)$;
$(ii)$ the $W(m;{\bm 1}_m)_0$-module $W(m;{\bm 1}_m)_{-1}$ is elementary; $(iii)$
$W(m;{\bm 1}_m)$ is the prolongation of $W(m;{\bm 1}_m)_-$.

We now consider the following dif\/ferential form
\[
\omega_K=dx_{2n+1}-\sum_{i=1}^nx_{i+n}dx_{i}.
\]
Def\/ine
\[
K(n)=\{ D\in W(2n+1):D\omega_K\in\mathfrak A(2n+1)\omega_K \}.
\]
(Here the action of $D$ on the dif\/ferential forms is extended from
its action $\mathfrak A(2n+1)$ by requiring that $D$ be derivation of
the exterior algebra satisfying $D(df)=d(Df)$, where
$df=\sum\frac{\partial f}{\partial x_i}dx_i$, $f\in\mathfrak A(m)$.)
We set $K(n)_p=W(2n+1;({\bm 1}_{2n},2))_p\cap K(n)$.
Then $K(n)=\bigoplus\limits_{p\geqq-2}K(n)_p$
is a transitive irreducible GLA
such that: $(i)$
$K(n)_0$ is isomorphic to $\csp(n,\mathbb C)$;
$(ii)$ the $K(n)_0$-module $K(n)_{-1}$ is elementary;
$(iii)$ $K(n)$ is the prolongation of $K(n)_-$
(cf.~\cite{kac68:1, mor88:1}).

From Proposition 2.2 of \cite{mt70:1}, we get
\begin{theorem}\label{thm6.1} Let $\gla g$ be a transitive GLA over $\mathbb C$
satisfying the following conditions:
\begin{enumerate}\itemsep=0pt
\renewcommand{\labelenumi}{$(\roman{enumi})$}
\item $\mathfrak g_-$ is an FGLA of the
$\mu$-th kind;
\item $\mathfrak g$ is infinite-dimensional;
\item The $\mathfrak g_0$-module $\mathfrak g_{-1}$ is irreducible;
\item $\mathfrak g$ is the prolongation of $\mathfrak g_-$.
\end{enumerate}
Then $\mu\leqq2$ and $\gla g$ is isomorphic to $W(m;{\bm 1}_m)$ or $K(n)$.
\end{theorem}

\section{Classif\/ication of the prolongations\\ of free
fundamental graded Lie algebras}\label{section7}

Let $\mathfrak m=\bigoplus\limits_{p<0}\mathfrak g_p$ be a free FGLA of type
$(n,\mu)$ over $\mathbb C$, and let $\mathfrak g(\mathfrak m)
=\bigoplus\limits_{p\in\mathbb Z}\mathfrak g(\mathfrak m)_p$ be the
prolongation of~$\mathfrak m$.
First of all, we assume that $\dim \mathfrak g(\mathfrak m)=\infty$.
By Theorem~\ref{thm6.1},
$\mathfrak g(\mathfrak m)$ is isomorphic to $K(m)$
as a GLA, where $n=2m$.
Since $K(m)_0$ is isomorphic to $\csp(m,\mathbb C)$ and since
$\mathfrak g(\mathfrak m)_0$ is isomorphic to $\gl(n,\mathbb C)$,
we see that $m=1$.
Therefore $\mathfrak g(\mathfrak m)$ is isomorphic to
$K(1)$ as a GLA.

Next we assume that $\dim \mathfrak g(\mathfrak m)<\infty$
and $\mathfrak g(\mathfrak m)_1\ne0$.
Since the $\mathfrak g(\mathfrak m)_0$-module $\mathfrak g(\mathfrak m)_{-1}$
is irreducible,
$\mathfrak g(\mathfrak m)$ is a f\/inite-dimensional simple GLA (see \cite{kn64:1,och70:1}).
By Theorem~\ref{thm5.2}, $\mathfrak g(\mathfrak m)$ is isomorphic to one of the following types:
\[
(B_l,\{\alpha_l\})\quad l\geqq3,\qquad (G_2,\{\alpha_1\}).
\]
Thus we get a proof of the following theorem:
\begin{theorem}\label{thm7.1}
Let $\mathfrak m=\bigoplus\limits_{p<0}\mathfrak g_p$ be a free FGLA of type
$(n,\mu)$ over $\mathbb C$, and let $\gla{g (\mathfrak m)}$ be
the prolongation of $\mathfrak m$.
Then one of the following cases occurs:
\begin{enumerate}\itemsep=0pt
\renewcommand{\labelenumi}{$(\alph{enumi})$}
\item $(n,\mu)\ne(n,2)$ $(n\geqq2)$, $(2,3)$. In this case,
$\mathfrak g(\mathfrak m)_1=\{0\}$.
\item
$(n,\mu)=
(n,2)$ $(n\geqq3)$, $(2,3)$.
In this case, $\dim\mathfrak g(\mathfrak m)<\infty$ and
$\mathfrak g(\mathfrak m)_1\ne\{0\}$. Furthermore
$\mathfrak g(\mathfrak m)$ is isomorphic to a finite-dimensional simple GLA of
type $(B_n,\{\alpha_n\})$ $(n\geqq3)$ or $(G_2,\{\alpha_1\})$ $(n=2)$.
\item $(n,\mu)=(2,2)$.
In this case,
$\dim\mathfrak g(\mathfrak m)=\infty$. Furthermore,
$\mathfrak g(\mathfrak m)$ is isomorphic to $K(1)$ as a GLA.
\end{enumerate}
\end{theorem}

\section{Free pseudo-product fundamental graded Lie algebras}\label{section8}

An FGLA $\mathfrak m=\bigoplus\limits_{p<0}\mathfrak g_p$ equipped with
nonzero subspaces $\mathfrak e$, $\mathfrak f$ of $\mathfrak g_{-1}$ is called
a pseudo-product FGLA if the following conditions hold:
\begin{enumerate}\itemsep=0pt
\renewcommand{\labelenumi}{$(\roman{enumi})$}
\item $\mathfrak g_{-1}=\mathfrak e\oplus\mathfrak f$;
\item $[\mathfrak e,\mathfrak e]=[\mathfrak f,\mathfrak f]=\{0\}$.
\end{enumerate}
The pair $(\mathfrak e,\mathfrak f)$ is called the pseudo-product structure of
the pseudo-product FGLA
$\mathfrak m=\bigoplus\limits_{p<0}\mathfrak g_p$.
We will also denote by the triplet $(\mathfrak m;\mathfrak e,\mathfrak f)$
the pseudo-product FGLA $\mathfrak m=\bigoplus\limits_{p<0}\mathfrak g_p$
with pseudo-product structure $(\mathfrak e,\mathfrak f)$.
Let $\mathfrak m=\bigoplus\limits_{p<0}\mathfrak g_p$ (resp.
$\mathfrak m'=\bigoplus\limits_{p<0}\mathfrak g'_p$) be a
pseudo-product FGLA with pseudo-product structure $(\mathfrak e,\mathfrak f)$
(resp. $(\mathfrak e',\mathfrak f')$).
We say that two pseudo-product FGLAs $(\mathfrak m;\mathfrak e,\mathfrak f)$
and $(\mathfrak m';\mathfrak e',\mathfrak f')$
are isomorphic if there exists a GLA isomorphism
$\varphi$ of $\mathfrak m$ onto $\mathfrak m'$ such that
$\varphi(\mathfrak e)=\mathfrak e'$ and $\varphi(\mathfrak f)=\mathfrak f'$.
\begin{proposition}\label{prop8.1} Let $\mathfrak m=\bigoplus\limits_{p<0}\mathfrak g_p$ be a
pseudo-product FGLA of the $\mu$-th kind with pseudo-product
structure $(\mathfrak e,\mathfrak f)$.
If $\mathfrak m$ is a free FGLA of type $(n,\mu)$,
then $n=2$.
\end{proposition}

\begin{proof} Let $(e_1,\dots,e_m)$ (resp.\ $(f_1,\dots,f_l)$) be a basis of
$\mathfrak e$ (resp.~$\mathfrak f$).
Since $[\mathfrak e,\mathfrak f]=\mathfrak g_{-2}$,
the space~$\mathfrak g_{-2}$ is generated by
$\{[e_i,f_j]:i=1,\dots,m,j=1,\dots,l\}$ as a vector space, so
$\dim \mathfrak g_{-2}\leqq ml$. On the other hand,
since $\mathfrak m$ is a free FGLA,
\[
\dim \mathfrak g_{-2}=\dim b(\mathfrak g_{-1},\mu)_{-2}=
\dim \Lambda^2(\mathfrak g_{-1})=\frac{(m+l)(m+l-1)}{2},
\]
so $ml\geqq \frac{(m+l)(m+l-1)}{2}$.
From this fact it follows that $m=l=1$.
\end{proof}

Let $\mathfrak m=\bigoplus\limits_{p<0}\mathfrak g_p$ be a pseudo-product FGLA
of the $\mu$-th kind with
pseudo-product structure $(\mathfrak e,\mathfrak f)$, where $\mu\geqq2$.
$\mathfrak m$ is called a free pseudo-product FGLA
of type $(m,n,\mu)$ if the following conditions hold:
\begin{enumerate}\itemsep=0pt
\renewcommand{\labelenumi}{$(\roman{enumi})$}
\item $\dim\mathfrak e=m$ and $\dim\mathfrak f=n$;
\item Let $\mathfrak m'=\bigoplus\limits_{p<0}\mathfrak g'_p$ be a pseudo-product FGLA of the $\mu$-th kind with
pseudo-product structure $(\mathfrak e',\mathfrak f')$
and let $\varphi$ be a surjective linear mapping of $\mathfrak g_{-1}$ onto
$\mathfrak g'_{-1}$ such that
$\varphi(\mathfrak e)\subset\mathfrak e'$ and
$\varphi(\mathfrak f)\subset\mathfrak f'$.
Then $\varphi$ can be extended uniquely to a GLA
epimorphism of $\mathfrak m$ onto $\mathfrak m'$.
\end{enumerate}
\begin{proposition}\label{prop8.2} Let $m$, $n$ and $\mu$ be positive integers such that $\mu\geqq2$.
\begin{enumerate}\itemsep=0pt
\item There exists a unique free pseudo-product FGLA
of type
$(m,n,\mu)$ up to isomorphism.
\item Let $\mathfrak m=\bigoplus\limits_{p<0}\mathfrak g_p$ be a free
pseudo-product FGLA of type $(m,n,\mu)$
with pseudo-product structure $(\mathfrak e,\mathfrak f)$.
We denote by $\Der(\mathfrak m;\mathfrak e,\mathfrak f)_0$ the Lie algebra of
all the derivations of $\mathfrak m$ preserving the gradation of $\mathfrak m$, $\mathfrak e$ and $\mathfrak f$.
Then the mapping $\Phi:\Der(\mathfrak m;\mathfrak e,\mathfrak f)_0\ni D\mapsto (D|\mathfrak e,D|\mathfrak f)\in\gl(\mathfrak e)\times\gl(\mathfrak f)$ is a Lie algebra isomorphism.
\end{enumerate}
\end{proposition}

\begin{proof} (1)
The uniqueness of a free pseudo-product FGLA of type $(m,n,\mu)$ follows from the def\/inition.
Let $V$ be an $(m+n)$-dimensional vector space and let
$\mathfrak e$, $\mathfrak f$ be subspaces of $V$
such that $V=\mathfrak e\oplus \mathfrak f$,
$\dim \mathfrak e=m$ and $\dim \mathfrak f=n$.
Let $\mathfrak a=\bigoplus\limits_{p<0}\mathfrak a_p$ be the graded ideal
of $b(V,\mu)$ generated by
$[\mathfrak e,\mathfrak e]+[\mathfrak f,\mathfrak f]$.
We set $\mathfrak m=b(V,\mu)/\mathfrak a$,
$\mathfrak g_p=b(V,\mu)_p/\mathfrak a_p$.
Clearly $\mathfrak m=\bigoplus\limits_{p<0}\mathfrak g_p$ is
a~pseudo-product FGLA.
We show that the factor algebra $\mathfrak m$
is a free pseudo-product FGLA of type $(m,n,\mu)$.
First we prove that $\mathfrak m$ is of the $\mu$-th kind.
Let $\mathfrak n=\bigoplus\limits_{p<0}\mathfrak g''_p$ be
a free FGLA of type $(2,\mu)$ and let~$\mathfrak e''$ and~$\mathfrak f''$
be one-dimensional subspaces of $\mathfrak g''_{-1}$ such that
$\mathfrak g''_{-1}=\mathfrak e''\oplus\mathfrak f''$.
Let~$\varphi_1$ be an injective linear mapping of $\mathfrak g''_{-1}$
into $V$ such that
$\varphi_1(\mathfrak e'')\subset \mathfrak e$ and
$\varphi_1(\mathfrak f'')\subset \mathfrak f$.
Let~$\varphi_2$ be a linear mapping of~$V$
into $\mathfrak g''_{-1}$ such that
$\varphi_2\circ\varphi_1=1_{\mathfrak g''_{-1}}$,
$\varphi_2(\mathfrak e)=\mathfrak e''$ and
$\varphi_2(\mathfrak f)=\mathfrak f''$.
There exists
a~homomorphism~$L(\varphi_1)$ (resp.~$L(\varphi_2)$) of
$\mathfrak n$ (resp.~$b(V,\mu)$) into
$b(V,\mu)$ (resp.~$\mathfrak n$)
such that
$L(\varphi_1)|\mathfrak g''_{-1}=\varphi_1$
(resp.\ $L(\varphi_2)|V=\varphi_2$).
Since $L(\varphi_2)([\mathfrak e,\mathfrak e]+[\mathfrak f,\mathfrak f])=\{0\}$,
$L(\varphi_2)$ induces a~homomorphism $\hat{L}(\varphi_2)$ of
$\mathfrak m$ into $\mathfrak n$ such that
$L(\varphi_2)=\hat{L}(\varphi_2)\circ \pi$,
where $\pi$ is the natural projection of $b(V,\mu)$ onto
$\mathfrak m$. Since
\[
1_{\mathfrak n}=L(\varphi_2)\circ L(\varphi_1)
=\hat{L}(\varphi_2)\circ \pi\circ L(\varphi_1),
\]
$\pi\circ L(\varphi_1)$ is a monomorphism of $\mathfrak n$ into $\mathfrak m$, so $\mathfrak g_{-\mu}\ne\{0\}$.
Thus $\mathfrak m$ is of the $\mu$-th kind.
Let $\mathfrak m'=\bigoplus\limits_{p<0}\mathfrak g'_p$ be a pseudo-product
FGLA of the $\mu$-th kind
with pseudo-product structure $(\mathfrak e',\mathfrak f')$
and let $\phi$ be a surjective linear mapping of $b(V,\mu)_{-1}$ onto
$\mathfrak g'_{-1}$ such that
$\phi(\mathfrak e)\subset\mathfrak e'$ and
$\phi(\mathfrak f)\subset\mathfrak f'$.
By the def\/inition of a free FGLA,
there exists a GLA epimorphism $L(\phi)$ of
$b(V,\mu)$ onto $\mathfrak m'$ such that
$L(\phi)|b(V,\mu)_{-1}=\phi$.
Since $L(\phi)([\mathfrak e,\mathfrak e]+[\mathfrak f,\mathfrak f])\subset [\mathfrak e',\mathfrak e']+[\mathfrak f',\mathfrak f']=\{0\}$,
we see that $L(\phi)(\mathfrak a)=\{0\}$, so
the epimorphism $L(\phi)$ induces a GLA epimorphism
$\hat{L}(\phi)$ of $\mathfrak m$ onto $\mathfrak m'$ such that
$\hat{L}(\phi)|\mathfrak g_{-1}=\phi$.

(2)~We may prove the fact that the mapping $\Phi$ is surjective.
Let $\phi$ be an endomorphism of~$\mathfrak g_{-1}$ such that
$\phi(\mathfrak e)\subset \mathfrak e$ and
$\phi(\mathfrak f)\subset \mathfrak f$.
By Proposition~\ref{prop2.1}~(2), there exists a $D\in\Der(b(V,\mu))_0$
such that $D|b(V,\mu)_{-1}=\phi$.
Since $D([\mathfrak e,\mathfrak e]+[\mathfrak f,\mathfrak f])
\subset [\mathfrak e,\mathfrak e]+[\mathfrak f,\mathfrak f]$,
$D$ induces a derivation~$\hat{D}$ of~$\mathfrak m$ such that
$\hat{D}|\mathfrak g_{-1}=\phi$.
\end{proof}

\begin{remark}\label{rem8.1}
Let $m$, $n$, $m'$, $n'$ and $\mu$ be positive integers with $\mu\geqq2$, and let
$\mathfrak m=\bigoplus\limits_{p<0}\mathfrak g_p$
(resp.\ $\mathfrak m'=\bigoplus\limits_{p<0}\mathfrak g'_p$) be a
free pseudo-product FGLA of type $(m,n,\mu)$ $($resp.~$(m',n',\mu))$
with pseudo-product structure $(\mathfrak e,\mathfrak f)$
(resp. $(\mathfrak e',\mathfrak f'))$.
Furthermore
let $\varphi$ be a linear mapping of $\mathfrak g_{-1}$ into
$\mathfrak g'_{-1}$ such that $\varphi(\mathfrak e)\subset \mathfrak e'$ and
$\varphi(\mathfrak f)\subset \mathfrak f'$.
\begin{enumerate}\itemsep=0pt
\item
From the proof of Proposition~\ref{prop8.2}, there exists a unique GLA
homomorphism $\hat{L}(\varphi)$ of $\mathfrak m$ into $\mathfrak m'$
such that $\hat{L}(\varphi)|\mathfrak g_{-1}=\varphi$.
If $\varphi$ is injective, then $\hat{L}(\varphi)$ is a monomorphism.
\item
Assume that $m=n=1$ and $\varphi$ is injective.
Then $\hat{L}(\varphi)(\mathfrak m)$ is a graded subalgebra of $\mathfrak m'$
isomorphic to a free FGLA of type $(2,\mu)$.
From this result,
the subalgebra of $\mathfrak m'$ generated by a nonzero element $X$ of
$\mathfrak e'$ and a nonzero element $Y$ of $\mathfrak f'$
is a free FGLA of type $(2,\mu)$.
\end{enumerate}
\end{remark}

Let $\mathfrak m=\bigoplus\limits_{p<0}\mathfrak g_p$ be a pseudo-product FGLA
of the $\mu$-th kind with pseudo-product structure~$(\mathfrak e,\mathfrak f)$.
We denote by $\mathfrak g_0$ the Lie algebra of
all the derivations of $\mathfrak m$ preserving the gradation of~$\mathfrak m$,~$\mathfrak e$ and~$\mathfrak f$:
\[
\mathfrak g_0=\{ D\in\Der(\mathfrak g)_0:D(\mathfrak e)\subset \mathfrak e,
D(\mathfrak f)\subset \mathfrak f \}.
\]
The prolongation
$\gla g$ of $(\mathfrak m,\mathfrak g_0)$ is called the prolongation of
$(\mathfrak m;\mathfrak e,\mathfrak f)$.

A transitive GLA $\gla g$ is called a pseudo-product GLA
if there are given nonzero subspaces $\mathfrak e$ and $\mathfrak f$ of $\mathfrak g_{-1}$ satisfying the following conditions:
\begin{enumerate}\itemsep=0pt
\renewcommand{\labelenumi}{$(\roman{enumi})$}
\item The negative part $\mathfrak g_{-}$ is
a pseudo-product
FGLA with pseudo-product structure
$(\mathfrak e,\mathfrak f)$;
\item $[\mathfrak g_0,\mathfrak e]\subset \mathfrak e$
and $[\mathfrak g_0,\mathfrak f]\subset \mathfrak f$.
\end{enumerate}
The pair $(\mathfrak e,\mathfrak f)$ is called the pseudo-product structure of
 the pseudo-product GLA
$\gla g$.
If the $\mathfrak g_0$-modules $\mathfrak e$ and $\mathfrak f$ are
irreducible,
then the pseudo-product GLA $\gla g$ is said to be
of irreducible type.

The following lemma is due to N.~Tanaka (cf.~\cite{tan89:01}).
Here we give a proof for the convenience of the readers.
\begin{lemma}\label{lemma8.1} Let $\gla g$ be a pseudo-product GLA of
depth $\mu$ with pseudo-product structure~$(\mathfrak e,\mathfrak f)$.
\begin{enumerate}\itemsep=0pt
\item If $\mathfrak g_-$ is non-degenerate,
then $\mathfrak g$ is finite-dimensional.
\item If $\gla g$ is of irreducible type and $\mu\geqq2$, then
$\mathfrak g$ is finite-dimensional.
\end{enumerate}
\end{lemma}

\begin{proof} (1) The proof is essentially due to the proof of
\cite[Corollary 3 to Theorem 11.1]{tan70:1}.
For $p\in \mathbb Z$, we set
$\mathfrak h_p=\{X\in\mathfrak g_p:[X,\mathfrak g_{\leqq -2}]=\{0\}\}$.
We def\/ine $I\in\mathfrak{gl}(\mathfrak g_{-1})$ as follows:
$I(x)=-\sqrt{-1}x$ for $x\in\mathfrak e$,
$I(x)=\sqrt{-1}x$ for $x\in\mathfrak f$.
Then $I^2=-1$,
$I([a,x])=[a,I(x)]$ and $[I(x),I(y)]=[x,y]$ for $a\in\mathfrak g_0$ and
$x,y\in\mathfrak g_{-1}$. We put $\langle x,y\rangle=[I(x),y]$ for
$x,y\in\mathfrak g_{-1}$. Then $\langle x,y\rangle=\langle y,x\rangle$,
and for $x\in\mathfrak g_{-1}$, $\langle x,\mathfrak g_{-1}\rangle=\{0\}$
implies $x=0$, since $\mathfrak g_-$ is non-degenerate. Also
$\langle[a,x],y\rangle+\langle x,[a,y]\rangle=0$ and $[[b,x],y]=[[b,y],x]$
for $a\in\mathfrak h_0$, $b\in\mathfrak h_1$ and $x,y\in\mathfrak g_{-1}$.
Then, for $b\in\mathfrak h_1$,
$x,y,z\in\mathfrak g_{-1}$, we have $\langle[[b,x],y],z\rangle=
-\langle y,[[b,x],z]\rangle=-\langle y,[[b,z],x]\rangle=
\langle [[b,z],y],x\rangle=\langle[[b,y],z],x\rangle
=-\langle z,[[b,y],x]\rangle=-\langle[[b,x],y],z\rangle$, so
$\langle[[b,x],y],z\rangle=0$. By transitivity of $\mathfrak g$,
$\mathfrak h_1=\{0\}$. Therefore by \cite[Corollary 1 to Theorem 11.1]{tan70:1}, $\mathfrak g$ is f\/inite-dimensional.

(2) We may assume that $\mathfrak g_1\ne\{0\}$.
By \cite[Lemma 2.4]{yat88:0}, the
$\mathfrak g_0$-modules $\mathfrak e$, $\mathfrak f$ are not isomorphic to each
other. We put
$\mathfrak d=\{X\in\mathfrak g_{-1}:[X,\mathfrak g_{-1}]=\{0\}\}$;
then $\mathfrak d$ is a $\mathfrak g_0$-submodule of $\mathfrak g_{-1}$.
Hence $\mathfrak d=\{0\}$, $\mathfrak d=\mathfrak e$, $\mathfrak d=\mathfrak f$ or
$\mathfrak d=\mathfrak g_{-1}$. If $\mathfrak d\ne\{0\}$, then $\mathfrak g_{-2}=[\mathfrak e,\mathfrak f]=\{0\}$, which is a contradiction. Thus $\mathfrak g_-$ is non-degenerate.
By (1), $\mathfrak g$ is f\/inite-dimensional.
\end{proof}

The prolongation of a pseudo-product FGLA becomes
a pseudo-product GLA.
By Proposition~\ref{prop8.2}~(2), the prolongation of a free pseudo-product
FGLA is a pseudo-product GLA of irreducible type.
By Lemma \ref{lemma8.1}~(2), the prolongation of a free pseudo-product FGLA is
f\/inite-dimensional.

\begin{proposition}\label{prop8.3} Let $\mathfrak m=\bigoplus\limits_{p<0}\mathfrak g_p$ be a
free pseudo-product FGLA of type $(m,n,\mu)$
with pseudo-product structure $(\mathfrak e,\mathfrak f)$
and let
$\gla g$ be the prolongation of $(\mathfrak m;\mathfrak e,\mathfrak f)$.
\begin{enumerate}\itemsep=0pt
\item
$\mathfrak g_0$ is isomorphic to $\gl(\mathfrak e)\oplus\gl(\mathfrak f)$ as
a Lie algebra.
\item $\mathfrak g_{-2}$ is isomorphic to $\mathfrak e\otimes \mathfrak f$
as a $\mathfrak g_0$-module. In particular,
$\dim\mathfrak g_{-2}=mn$.
\item $\mathfrak g_{-3}$ is isomorphic to
$S^2(\mathfrak e)\otimes
\mathfrak f\oplus S^2(\mathfrak f)\otimes \mathfrak e$
as a $\mathfrak g_0$-module.
In particular,
$\dim\mathfrak g_{-3}=\frac{mn(m+n+2)}2$.
\end{enumerate}
\end{proposition}

\begin{proof} (1) This follows from Proposition \ref{prop8.2} (2).
\par
(2) Let $\mathfrak a=\bigoplus\limits_{p<0}\mathfrak a_p$ be the graded ideal
of $b(\mathfrak g_{-1},\mu)$ generated by
$[\mathfrak e,\mathfrak e]+[\mathfrak f,\mathfrak f]$.
By the construction of $b(\mathfrak g_{-1},\mu)_{-2}$,
$\mathfrak a_{-2}$ is isomorphic to $\Lambda^2(\mathfrak e)\oplus
\Lambda^2(\mathfrak f)$, so
$\mathfrak g_{-2}=b(\mathfrak g_{-1},\mu)_{-2}/\mathfrak a_{-2}$
is isomorphic to $\mathfrak e\otimes \mathfrak f$.

(3) By the construction of $b(\mathfrak g_{-1},\mu)_{-3}$,
$b(\mathfrak g_{-1},\mu)_{-3}$ is isomorphic to
\[
(\mathfrak e\oplus \mathfrak f)\otimes \Lambda^2(\mathfrak e\oplus \mathfrak f)/\Lambda^3(\mathfrak e\oplus \mathfrak f)\cong
(\mathfrak e\otimes \mathfrak e\otimes \mathfrak f)
\oplus
(\mathfrak e\otimes \mathfrak f\otimes \mathfrak f).
\]
Moreover,
$\mathfrak a_{-3}$ is isomorphic to
\[
(\mathfrak e\oplus \mathfrak f)\otimes \Lambda^2(\mathfrak e)
\oplus(\mathfrak e\oplus \mathfrak f)\otimes \Lambda^2(\mathfrak f)/
\Lambda^3(\mathfrak e\oplus \mathfrak f)
\cong \mathfrak e\otimes\Lambda^2(\mathfrak e)\oplus
\mathfrak f\otimes\Lambda^2(\mathfrak f).
\]
Hence
$\mathfrak g_{-3}=b(\mathfrak g_{-1},\mu)_{-3}/\mathfrak a_{-3}$
is isomorphic to
\[
(\mathfrak e\otimes \mathfrak e\otimes \mathfrak f)/
\Lambda^2(\mathfrak e)\otimes\mathfrak f\oplus
(\mathfrak e\otimes \mathfrak f\otimes \mathfrak f)/
\mathfrak e\otimes\Lambda^2(\mathfrak f)
\cong
S^2(\mathfrak e)\otimes
\mathfrak f\oplus S^2(\mathfrak f)\otimes \mathfrak e.
\]
This completes the proof.
\end{proof}

\begin{proposition}\label{prop8.4}
Let $\mathfrak m=\bigoplus\limits_{p<0}\mathfrak g_p$ be a
pseudo-product FGLA of the $\mu$-th kind with
pseudo-product structure $(\mathfrak e,\mathfrak f)$, where $\mu\geqq2$.
We denote by $\mathfrak c$
the centralizer of $\mathfrak g_{-2}$ in $\mathfrak g_{-1}$.
Let $\gla g$ be the prolongation of $(\mathfrak m;\mathfrak e,\mathfrak f)$.
Assume that $\mathfrak g_0$ is isomorphic to $\gl(\mathfrak e)\oplus\gl(\mathfrak f)$ as a Lie algebra.
\begin{enumerate}\itemsep=0pt
\item If $\mu=2$, then $\mathfrak m=\bigoplus\limits_{p<0}\mathfrak g_p$
be a free
pseudo-product FGLA.
\item If $\mu\geqq3$ and $\mathfrak c\ne\{0\}$,
then $(\mathfrak m;\mathfrak e,\mathfrak f)$ is not a free pseudo-product
FGLA.
\item If $\mu=3$ and $\mathfrak c=\{0\}$,
then $(\mathfrak m;\mathfrak e,\mathfrak f)$ is a free pseudo-product
FGLA.
\end{enumerate}
\end{proposition}

\begin{proof}
Let
$\check{\mathfrak m}=\bigoplus_{p<0}\limits\check{\mathfrak g}_p$
be the free pseudo-product FGLA of type $(m,n,\mu)$
with pseudo-product structure
$(\check{\mathfrak e},\check{\mathfrak f})$ such that
$\check{\mathfrak g}_{-1}=\mathfrak g_{-1}$,
$\check{\mathfrak e}=\mathfrak e$ and $\check{\mathfrak f}=\mathfrak f$.
Let $\check{\mathfrak g}=\bigoplus_{p\in\mathbb Z}\limits\check{\mathfrak g}_p$  be the prolongation of $(\check{\mathfrak m};\check{\mathfrak e},\check{\mathfrak f})$.
There exists a GLA epimorphism $\varphi$
of $\check{\mathfrak m}$ onto $\mathfrak m$ such that the restriction
$\varphi|\check{\mathfrak g}_{-1}$ is the identity mapping.
Since the mapping
$\check{\mathfrak g}_0\ni D\mapsto
(D|\mathfrak e,D|\mathfrak f)\in\gl(\mathfrak e)\times\gl(\mathfrak f)$
is an isomorphism,
$\varphi$ can be extended to be a homomorphism $\check{\varphi}$
of $\bigoplus\limits_{p\leqq0}\check{\mathfrak g}_{p}$ onto
$\bigoplus\limits_{p\leqq0}\mathfrak g_p$.
Let $\mathfrak a$ be the kernel of $\check{\varphi}$;
then $\mathfrak a$ is a graded ideal of
$\bigoplus\limits_{p\leqq0}\check{\mathfrak g}_{p}$.
We set $\mathfrak a_p
=\mathfrak a\cap \check{\mathfrak g}_{p}$; then
$\mathfrak a=\bigoplus\limits_{p\leqq 0}\mathfrak a_p$.
Since the restriction of $\check{\varphi}$ to
$\check{\mathfrak g}_{-1}\oplus \check{\mathfrak g}_0$ is injective,
$\mathfrak a_p=\{0\}$ for $p\geqq-1$. Also
each $\mathfrak a_p$ is a $\check{\mathfrak g}_0$-submodule of
$\check{\mathfrak g}_p$.
Since the $\check{\mathfrak g}_0$-module
$\check{\mathfrak g}_{-2}$ is irreducible (Proposition~\ref{prop8.3}~(2)),
$\varphi|\mathfrak g_{-2}$ is injective.
If $\mu=2$, then $\varphi$ is an isomorphism. This proves the assertion~(1).
Now we assume that $\mu\geqq3$. Then
\[
\check{\mathfrak g}_{-3}=[[\mathfrak e,\mathfrak f],\mathfrak f]\oplus
[[\mathfrak e,\mathfrak f],\mathfrak e].
\]
Since $\check{\mathfrak g}_0$-modules
$[[\mathfrak e,\mathfrak f],\mathfrak f]$ and
$[[\mathfrak e,\mathfrak f],\mathfrak e]$ are irreducible and not isomorphic
to each other (Proposition \ref{prop8.3} (3)),
one of the following cases occurs:
$(i)$~$\mathfrak a_{-3}=[[\mathfrak e,\mathfrak f],\mathfrak f]$;
$(ii)$~$\mathfrak a_{-3}=[[\mathfrak e,\mathfrak f],\mathfrak e]$;
$(iii)$~$\mathfrak a_{-3}=\{0\}$.
If $\mathfrak a_{-3}=[[\mathfrak e,\mathfrak f],\mathfrak f]$ (resp.\
$\mathfrak a_{-3}=[[\mathfrak e,\mathfrak f],\mathfrak e]$), then
$\mathfrak c=\mathfrak f$
(resp. $\mathfrak c=\mathfrak e$).
Also
since $\mathfrak g_0$-modules
$\mathfrak e$,
$\mathfrak f$ are irreducible and not isomorphic to each other,
one of the following cases occurs:
$(i)$~$\mathfrak c=\mathfrak e$;
$(ii)$~$\mathfrak c=\mathfrak f$;
$(iii)$~$\mathfrak c=\{0\}$.
If $\mathfrak c=\mathfrak e$
(resp. $\mathfrak c=\mathfrak f$),
then $\mathfrak a_{-3}=[[\mathfrak e,\mathfrak f],\mathfrak e]$
(resp. $\mathfrak a_{-3}=[[\mathfrak e,\mathfrak f],\mathfrak f]$).
In this case, $\varphi$ is not injective.
Hence $(\mathfrak m;\mathfrak e,\mathfrak f)$ is not free.
If $\mathfrak c=\{0\}$,
then $\mathfrak a_{-3}=\{0\}$.
Hence $\varphi|\check{\mathfrak g}_{-3}$ is an isomorphism.
In particular, if $\mu=3$, then
$(\mathfrak m;\mathfrak e,\mathfrak f)$ is free.
\end{proof}

\begin{example}
Let $V$ and $W$ be f\/inite-dimensional vector spaces and $k\geqq1$.
We set
\begin{gather*}
\mathfrak C^{k}(V,W) =\bigoplus\limits_{p=-k-1}^{-1}\mathfrak C^{k}(V,W)_p, \\
\mathfrak C^{k}(V,W)_p=W\otimes S^{k+p+1}(V^*), \qquad
 -k-1\leqq p\leqq -2, \\
\mathfrak C^{k}(V,W)_{-1} =V\oplus (W\otimes S^{k}(V^*)).
\end{gather*}
The bracket operation of $\mathfrak C^{k}(V,W)$ is def\/ined as follows:
\begin{gather*}
  [W,V]=\{0\},\qquad [V,V]=\{0\}, \qquad [W\otimes S^{r}(V^*),W\otimes S^{s}(V^*)]=\{0\}, \\
[w\otimes s_r,v]=w\otimes (v\lrcorner\,s_r)\qquad \text{for}\  v\in V,\ w\in W,
\ s_r\in S^r(V^*).
\end{gather*}
Equipped with this bracket operation, $\mathfrak C^{k}(V,W)$ becomes
a pseudo-product FGLA of the $(k+1)$-th kind
with pseudo-product structure $(V,W\otimes S^{k}(V^*))$,
which is called
{\it the contact algebra of order $k$ of bidegree} $(n,m)$, where
$n=\dim V$ and $m=\dim W$ (cf.~\cite[p.~133]{yam82:1}).
We assume that $\mathfrak C^{k}(V,W)$ is a free pseudo-product FGLA.
Since
\begin{gather*}
\dim \mathfrak C^{k}(V,W)_{-2}=m\binom{n+k-2}{k-1}, \qquad
\dim V\dim (W\otimes S^{k}(V^*))= nm\binom{n+k-1}{k},
\end{gather*}
we get $n=1$.
Since $W\otimes S^{k}(V^*)$ is contained in the centralizer of
$\mathfrak C^{k}(V,W)_{-2}$
in $\mathfrak C^{k}(V,W)_{-1}$, we get $k=1$.
Thus we obtain that $\mathfrak C^{k}(V,W)$ is a free pseudo-product FGLA
if and only if $k=1$, $n=1$.
\end{example}

\begin{example}\label{ex8.2}
Let $\gla g$ be a f\/inite-dimensional simple GLA of type
$(A_{m+n},\{\alpha_m,\alpha_{m+1}\})$.
We set $\mathfrak e=\mathfrak g_{-1}^{(m)}$, $\mathfrak f=\mathfrak g_{-1}^{(m+1)}$. Then $(\mathfrak g_-;\mathfrak e,\mathfrak f)$
is a pseudo-product FGLA.
Since $\dim \mathfrak e=m$, $\dim\mathfrak f=n$ and
$\dim\mathfrak g_{-2}=mn$, the pseudo-product FGLA
$(\mathfrak g_-;\mathfrak e,\mathfrak f)$ is a free pseudo-product FGLA
of type $(m,n,2)$ (Proposition 8.3 (2)). Also
$\gla g$ is the prolongation of $\mathfrak g_-$ except for the following
cases (see \cite{yam93:1}):
\begin{enumerate}\itemsep=0pt
\item $m=n=1$. In this case,
the prolongation of $\mathfrak g_-$ is isomorphic to $K(1)$.
\item $m=1$ or $n=1$ and $l=\max\{m,n\}\geqq2$. In this case,
the prolongation of $\mathfrak g_-$ is isomorphic to $W(l+1;\bm s)$,
where $\bm s=(1,2,\dots,2)$.
\end{enumerate}
\end{example}

\begin{example}
Let $V$ and $W$ be f\/inite-dimensional vector spaces
such that $\dim V=m\geqq1$ and $\dim W=n\geqq1$.
We set
\begin{gather*}
\mathfrak g_{-1}  =V\oplus W, \qquad \mathfrak g_{-2}=V\otimes W, \\
\mathfrak g_{-3}  = V\otimes S^2(W)\oplus S^2(V)\otimes W,
\qquad \mathfrak m=\mathfrak g_{-1}\oplus \mathfrak g_{-2}\oplus
\mathfrak g_{-3}.
\end{gather*}
The bracket operation of $\mathfrak m$ is def\/ined as follows:
 \begin{gather*}
  [\mathfrak g_{-3},\mathfrak g_{-1}\oplus \mathfrak g_{-2}]
=[\mathfrak g_{-2},\mathfrak g_{-2}]=\{0\},  \qquad [V,V]=[W,W]=\{0\},\\
[v,w]=-[w,v]=v\otimes w, \qquad
  [v,v'\otimes w]=-[v'\otimes w,v]=v\circledcirc v'\otimes w,
\\ [v\otimes w,w']=-[w',v\otimes w]=v\otimes w\circledcirc w',
\end{gather*}
where $v,v'\in V$ and $w,w'\in W$.
Equipped with this bracket operation, $\mathfrak m$ becomes
a free pseudo-product FGLA of type $(m,n,3)$ with pseudo-product structure
$(V,W)$ (Proposition~\ref{prop8.3}).
\end{example}

\begin{theorem}\label{thm8.1}
Let $\mathfrak m=\bigoplus\limits_{p<0}\mathfrak g_p$ be a free
pseudo-product FGLA of type $(m,n,\mu)$ with
pseudo-product structure $(\mathfrak e,\mathfrak f)$ over $\mathbb C$.
Furthermore let
$\gla g$
$($resp.\
$\mathfrak g(\mathfrak m)
=\bigoplus\limits_{p\in\mathbb Z}\mathfrak g(\mathfrak m)_p)$
be the prolongation of $(\mathfrak m;\mathfrak e,\mathfrak f)$
$($resp.~$\mathfrak m)$.
\begin{enumerate}\itemsep=0pt
\item Assume that $\dim \mathfrak g(\mathfrak m)=\infty$.
Then $m=1$ or $n=1$, and $\mu=2$.
Furthermore $\gla g$ is isomorphic to
a finite-dimensional simple GLA of type $(A_{l+1},\{\alpha_1,\alpha_2\})$,
where $l=\max\{m,n\}$.
If $l=1$, then $\mathfrak g(\mathfrak m)$
is isomorphic to $K(1)$.
If $l\geqq2$, then $\mathfrak g(\mathfrak m)$
is isomorphic to $W(l+1;\bm s)$,
where $\bm s=(1,2,\dots,2)$.
\item
If $\mathfrak g_{1}\!\ne\!\{0\}$, then $\gla g$
is a finite-dimensional simple GLA of type
$(A_{m{+}n},\{\alpha_m,\alpha_{m{+}1}\}).\!$
\end{enumerate}
\end{theorem}

\begin{proof}
(1) For $p\!\geqq\!{-}1$, we put
$\mathfrak h_p
=\{X\in\mathfrak g(\mathfrak m)_p:[X,\mathfrak g_{\leqq -2}]=\{0\}\}$.
Assume that \mbox{$\dim \mathfrak g(\mathfrak m)=\infty$} and $\mu\geqq3$.
By Proposition \ref{prop8.4} (2),
$\mathfrak h_{-1}=\{0\}$ .
Since $[\mathfrak h_0,\mathfrak g_{-1}]\subset \mathfrak h_{-1}=\{0\}$,
we see that $\mathfrak h_0=\{0\}$.
By \cite[Corollary 1 to Theorem 11.1]{tan70:1}, we obtain that
$\dim\mathfrak g(\mathfrak m)<\infty$,
which is a~contradiction. Thus we see that $\mu=2$ if $\dim \mathfrak g(\mathfrak m)=\infty$. The remaining assertion follows from Example~\ref{ex8.2}.

(2)
Assume that $\mathfrak g_1\ne\{0\}$ and $\mu\geqq3$.
By transitivity of $\mathfrak g$,
$[\mathfrak g_1,\mathfrak e]\ne\{0\}$ or $[\mathfrak g_1,\mathfrak f]\ne\{0\}$.
We may assume that $[\mathfrak g_1,\mathfrak e]\ne\{0\}$.
Then there exists an irreducible component $\mathfrak g'_1$ of
the $\mathfrak g_0$-module~$\mathfrak g_1$
such that $[\mathfrak g'_1,\mathfrak e]\ne\{0\}$ and
$[\mathfrak g'_1,\mathfrak f]=\{0\}$.
The subalgebra $\mathfrak e\oplus[\mathfrak e,\mathfrak g'_1]\oplus
\mathfrak g'_1$ is a simple GLA of the f\/irst kind.
Since $\mathfrak g_0$ is isomorphic to $\mathfrak{gl}(\mathfrak e)\oplus
\mathfrak{gl}(\mathfrak f)$,
$\mathfrak e\oplus[\mathfrak e,\mathfrak g'_1]\oplus \mathfrak g'_1$ is of type $(A_m,\{\alpha_1\})$.
Let~$D$ be a nonzero element of $\mathfrak g'_1$.
There exist $\lambda\in \mathfrak e^*$ and $\eta\in\mathfrak f^*$ such that
\[
[[D,Z],U]=\lambda(U)Z+\lambda(Z)U, \qquad [[D,Z],W]=\eta(Z)W,
\]
where $Z,U\in \mathfrak e$ and $W\in \mathfrak f$
(cf.~\cite[p.~4]{tan57:1}).
Let $X$ (resp.~$Y$) be a nonzero element of $\mathfrak e$
(resp. $\mathfrak f$).
Since the subalgebra generated by $X,Y$ is a free FGLA of type $(2,\mu)$
(Remark~\ref{rem8.1}~(2)),
\begin{gather*}
\ad(X)^\mu(Y)=0, \qquad \ad(X)^{\mu-1}(Y)\ne 0, \\
\ad(Y)\ad(X)^{\mu-1}(Y)=0, \qquad \ad(Y)\ad(X)^{\mu-2}(Y)\ne 0
\end{gather*}
(Lemma \ref{lem2.1}).
By induction on $\mu$, we see that
\begin{gather*}
0  = \ad(D)\ad(X)^\mu(Y)=(\mu(\mu-1)\lambda(X)+\mu\eta(X))\ad(X)^{\mu-1}(Y),
\\
0  = \ad(D)\ad(Y)\ad(X)^{\mu-1}(Y)\\
\hphantom{0}{} =((\mu-1)(\mu-2)\lambda(X)
+(\mu-1)\eta(X))\ad(Y)\ad(X)^{\mu-2}(Y).
\end{gather*}
Since
\[
\det\begin{bmatrix}
\mu(\mu-1) & \mu \\
(\mu-1)(\mu-2) & \mu-1
\end{bmatrix}=
\mu(\mu-1)\ne0,
\]
we see that $\lambda(X)=\eta(X)=0$, which is a contradiction.
Thus we obtain that $\mu=2$ if \mbox{$\dim \mathfrak g_1\ne\{0\}$}.
From Example~\ref{ex8.2}, it follows that
$\gla g$ is a simple GLA of type $(A_{m+n}$, $\{\alpha_m,\alpha_{m+1}\})$
if $\dim \mathfrak g_1\ne\{0\}$.
\end{proof}

\section{Automorphism groups of the prolongations\\ of free pseudo-product
fundamental graded Lie algebras}\label{section9}

For a GLA $\gla g$
we denote by $\Aut(\mathfrak g)_0$ the group of all the automorphisms of
$\mathfrak g$ preserving the gradation of $\mathfrak g$:
\[
\Aut(\mathfrak g)_0=\{\varphi\in\Aut(\mathfrak g):\varphi(\mathfrak g_p)=
\mathfrak g_p \ \text{for all}\ p\in\mathbb Z\}.
\]

\begin{proposition}\label{prop9.1} Let $\mathfrak m=\bigoplus\limits_{p<0}\mathfrak g_p$ be
an FGLA and let $\mathfrak g(\mathfrak m)
=\bigoplus\limits_{p\in\mathbb Z}\mathfrak g(\mathfrak m)_p$
be the prolongation of~$\mathfrak m$.
The mapping
$\Phi:\Aut(\mathfrak g(\mathfrak m))_0\ni\phi\mapsto
\phi|\mathfrak m\in\Aut(\mathfrak m)_0$ is an isomorphism.
\end{proposition}

\begin{proof}
It is clear that $\Phi$ is a group homomorphism.
We prove that $\Phi$ is injective.
Let $\phi$ be an element of $\Ker\Phi$.
Assume that $\phi(X)=X$ for all $X\in \mathfrak g(\mathfrak m)_p$ $(p<k)$.
For $X\in\mathfrak g(\mathfrak m)_k$, $Y\in\mathfrak g_{-1}$,
\[
[\phi(X)-X,Y]=\phi([X,Y])-[X,Y].
\]
Since $[X,Y]\in\mathfrak g(\mathfrak m)_{k-1}$, we have $[\phi(X)-X,Y]=0$.
By transitivity, $\phi(X)=X$.
By induction, we have proved $\phi$ to be the identity mapping.
Hence $\Phi$ is a monomorphism.

We prove that $\Phi$ is surjective.
Let $\varphi\in\Aut(\mathfrak m)_0$.
We construct the mapping
$\varphi_p:\mathfrak g(\mathfrak m)_p\to\mathfrak g(\mathfrak m)_p$
inductively as follows:
First for $X\in\mathfrak g(\mathfrak m)_0$, we set
$\varphi_0(X)=\varphi X\varphi^{-1}$.
Then for $Y,Z\in \mathfrak m$
\[
\varphi_0(X)([Y,Z])=[\varphi(X(\varphi^{-1}(Y))),Z]
+[Y,\varphi(X(\varphi^{-1}(Z)))],
\]
so $\varphi_0(X)\in\mathfrak g(\mathfrak m)_0$.
Furthermore we can prove easily that
$[\varphi_0(X),\varphi_p(Y)]=\varphi_p([X,Y])$ for
$X\in\mathfrak g_0$ and $Y\in \mathfrak g_p$ $(p\leqq 0)$.
Here for $p<0$ we set $\varphi_p=\varphi|\mathfrak g(\mathfrak m)_p$.
Assume that we have def\/ined linear isomorphisms $\varphi_p$ of
$\mathfrak g(\mathfrak m)_p$ onto itself $(0\leqq p<k)$ such that
\[
\varphi_{r+s}([X,Y])=[\varphi_r(X),\varphi_s(Y)]
\]
for $X\in\mathfrak g(\mathfrak m)_r$, $Y\in\mathfrak g(\mathfrak m)_s$
$(r+s<k$, $r<k$, $s<k)$.
For $X\in\mathfrak g(\mathfrak m)_k$ we def\/ine $\varphi_k(X)\in
\Hom(\mathfrak m,\bigoplus\limits_{p\leqq k-1}\mathfrak g(\mathfrak m)_p)_k$
as follows:
\[
\varphi_k(X)(Y)=\varphi_{k+s}([X,\varphi^{-1}(Y)]), \qquad Y\in\mathfrak g_s, \ s<0.
\]
For $Y\in\mathfrak g_s$, $Z\in\mathfrak g_t$ $(s,t<0)$,
\begin{gather*}
\varphi_k(X)([Y,Z]) =\varphi_{k+t+s}([X,\varphi^{-1}([Y,Z]]) \\
\hphantom{\varphi_k(X)([Y,Z])}{}
=\varphi_{k+s+t}([[X,\varphi^{-1}(Y)],\varphi^{-1}(Z)]
+[\varphi^{-1}(Y),[X,\varphi^{-1}(Z)]]) \\
\hphantom{\varphi_k(X)([Y,Z])}{}
=[\varphi_{k+s}([X,\varphi^{-1}(Y)]),Z]
+[Y,\varphi_{k+t}([X,\varphi^{-1}(Z)])] \\
\hphantom{\varphi_k(X)([Y,Z])}{}
=[\varphi_k(X)(Y),Z]+[Y,\varphi_k(X)(Z)],
\end{gather*}
so $\varphi_k(X)\in\mathfrak g(\mathfrak m)_k$.
Next we prove that for $X\in\mathfrak g_p$, $Y\in\mathfrak g_q$ $(p+q=k$,
$0\leqq p\leqq k$, $0\leqq q\leqq k)$,
\[
\varphi_k([X,Y])=[\varphi_p(X),\varphi_q(Y)].
\]
For $Z\in\mathfrak g_s$ $(s<0)$,
\begin{gather*}
 [[\varphi_p(X),\varphi_q(Y)],Z]
=[\varphi_p(X),[\varphi_q(Y),Z]]-[\varphi_q(Y),[\varphi_p(X),Z]] \\
\hphantom{[[\varphi_p(X),\varphi_q(Y)],Z]}{}
=\varphi_{p+q+s}([X,[Y,\varphi^{-1}(Z)]]-[Y,[X,\varphi^{-1}(Z)]]) \\
\hphantom{[[\varphi_p(X),\varphi_q(Y)],Z]}{}
=\varphi_{p+q+s}([[X,Y],\varphi^{-1}(Z)])
=[\varphi_k([X,Y]),Z].
\end{gather*}
By transitivity, we see that
$\varphi_k([X,Y])=[\varphi_p(X),\varphi_q(Y)]$.
We def\/ine a mapping $\check{\varphi}$ of $\mathfrak g(\mathfrak m)$ into
itself as follows:
\[\check{\varphi}(X)=
\begin{cases}
\varphi(X), & X\in\mathfrak m, \\
\varphi_k(X), & k\geqq0, \ X\in\mathfrak g(\mathfrak m)_k.
\end{cases}
\]
From the above results and the def\/inition of $\varphi_k$ $(k\geqq0)$,
we see that $\check{\varphi}$ is a GLA homomorphism.

Assume that $\varphi_{k-1}$ $(k\geqq0)$ is a linear isomorphism.
For $X\in \mathfrak g(\mathfrak m)_k$, if $\varphi_k(X)=0$, then
$0=[\varphi_k(X),Y]=\varphi_{k-1}([X,\varphi^{-1}(Y)])$ for all
$Y\in\mathfrak g_{-1}$. By transitivity, we see that $X=0$,
so $\varphi_k$ is a linear isomorphism.
Therefore $\check{\varphi}$ is an automorphism of $\mathfrak g(\mathfrak m)$.
\end{proof}

\begin{theorem}\label{thm9.1} Let $\mathfrak m=\bigoplus\limits_{p<0}\mathfrak g_p$ be a
free FGLA over $\mathbb C$, and let
$\mathfrak g(\mathfrak m)
=\bigoplus\limits_{p\in\mathbb Z}\mathfrak g(\mathfrak m)_p$
be the prolongation of $\mathfrak m$.
The mapping $\Phi:\Aut(\mathfrak g(\mathfrak m))_0\ni\phi\mapsto \phi|\mathfrak g_{-1}\in GL(\mathfrak g_{-1})$ is an isomorphism.
\end{theorem}

\begin{proof}
We may assume that $\mathfrak m$ is a universal FGLA $b(\mathfrak g_{-1},\mu)$ of the $\mu$-th kind.
By Corollary~1 to Proposition~3.2 of~\cite{tan70:1},
the mapping $\Aut(\mathfrak m)_0\ni a\mapsto a|\mathfrak g_{-1}\in
GL(\mathfrak g_{-1})$ is an isomorphism.
By Proposition~\ref{prop9.1}, we see that the mapping
$\Phi:\Aut(\mathfrak g(\mathfrak m))_0\ni\phi\mapsto \phi|\mathfrak g_{-1}\in GL(\mathfrak g_{-1})$ is an isomorphism.
\end{proof}

For a pseudo-product GLA $\gla g$ with pseudo-product structure
$(\mathfrak e,\mathfrak f)$,
we denote by $\Aut(\mathfrak g;\mathfrak e,\mathfrak f)_0$ the group of all the automorphisms of
$\mathfrak g$ preserving the gradation of $\mathfrak g$, $\mathfrak e$
and $\mathfrak f$:
\[
\Aut(\mathfrak g;\mathfrak e,\mathfrak f)_0
=\{\varphi\in\Aut(\mathfrak g)_0:
\varphi(\mathfrak e)=\mathfrak e,\varphi(\mathfrak f)=\mathfrak f
\}.
\]

\begin{theorem}\label{thm9.2} Let $\mathfrak m=\bigoplus\limits_{p<0}\mathfrak g_p$ be a
free pseudo-product FGLA of type $(m,n,\mu)$ with pseudo-product structure
$(\mathfrak e,\mathfrak f)$ over $\mathbb C$, and let
$\mathfrak g
=\bigoplus\limits_{p\in\mathbb Z}\mathfrak g_p$
be the prolongation of $(\mathfrak m;\mathfrak e,\mathfrak f)$.
The mapping $\Phi:\Aut(\mathfrak g;\mathfrak e,\mathfrak f)_0\ni\phi\mapsto
(\phi|\mathfrak e,\phi|\mathfrak f)\in GL(\mathfrak e)\times GL(\mathfrak f)$
is an isomorphism.
Furthermore if $\dim\mathfrak e\ne\dim \mathfrak f$,
then $\Aut(\mathfrak g;\mathfrak e,\mathfrak f)_0=\Aut(\mathfrak g)_0$.
\end{theorem}

\begin{proof}
Clearly $\Phi$ is a monomorphism.
We show that $\Phi$ is surjective.
Let $(\phi_1,\phi_2)$ be an element of
$GL(\mathfrak e)\times GL(\mathfrak f)$.
We set $\phi=\phi_1\oplus\phi_2\in GL(\mathfrak g_{-1})$.
By Corollary~1 to Proposition~3.2 of \cite{tan70:1},
there exists an element
$\varphi_1\in \Aut(b(\mathfrak g_{-1},\mu))_0$ such that
$\varphi_1|\mathfrak g_{-1}=\phi$.
Since $\varphi_1([\mathfrak e,\mathfrak e]+[\mathfrak f,\mathfrak f])=
[\mathfrak e,\mathfrak e]+[\mathfrak f,\mathfrak f]$,
$\varphi_1$ induces an element
$\varphi_2\in \Aut(\mathfrak m;\mathfrak e,\mathfrak f)_0$
such that $\varphi_2|\mathfrak g_{-1}=\phi$.
By Proposition \ref{prop9.1}, there exists
$\varphi_3\in \Aut(\mathfrak g(\mathfrak m))_0$
such that $\varphi_3|\mathfrak m=\varphi_2$.
We prove that $\varphi_3(\mathfrak g)=\mathfrak g$.
For $X_0\in\mathfrak g_0$ and $Y\in\mathfrak e$,
we see that $[\varphi_3(X_0),Y]=\varphi_3([X_0,\varphi_3^{-1}(Y)])\in
\varphi_3(\mathfrak e)=\mathfrak e$, so
$\varphi_3(X_0)(\mathfrak e)\subset \mathfrak e$.
Similarly we get $\varphi_3(X_0)(\mathfrak f)\subset \mathfrak f$.
Thus we obtain that $\varphi_3(\mathfrak g_0)=\mathfrak g_0$.
Now we assume that $\varphi_i(\mathfrak g_i)=\mathfrak g_i$
for $0\leqq i\leqq k$. Then
for $X_{k+1}\in\mathfrak g_{k+1}$ and $Y\in\mathfrak g_{p}$ $(p<0)$,
we see that $[\varphi_3(X_{k+1}),Y]=\varphi_3([X_{k+1},\varphi_3^{-1}(Y)])\in
\varphi_3(\mathfrak g_{p+k+1})=\mathfrak g_{p+k+1}$, so
$\varphi_3(\mathfrak g_{k+1})\subset \mathfrak g_{k+1}$.
Hence $\varphi_3(\mathfrak g)=\mathfrak g$ and $\Phi$ is surjective.
Now we assume that $\dim\mathfrak e\ne\dim \mathfrak f$.
Let $\varphi\in\Aut(\mathfrak g)_0$.
Since $\mathfrak g_0$-modules $\mathfrak e$ and $\mathfrak f$ are not
isomorphic to each other, we see that
$(i)$ $\varphi(\mathfrak e)=\mathfrak e$, $\varphi(\mathfrak f)=\mathfrak f$ or
$(ii)$ $\varphi(\mathfrak e)=\mathfrak f$, $\varphi(\mathfrak f)=\mathfrak e$.
According to the assumption $\dim\mathfrak e\ne\dim\mathfrak f$,
we get $\varphi(\mathfrak e)=\mathfrak e$, $\varphi(\mathfrak f)=\mathfrak f$,
so $\varphi\in \Aut(\mathfrak g;\mathfrak e,\mathfrak f)_0$.
\end{proof}

\pdfbookmark[1]{References}{ref}
\LastPageEnding

\end{document}